\documentclass{article}

\usepackage{graphicx}
\usepackage[utf8]{inputenc}
\usepackage{textcomp}

\usepackage{amsmath}
\usepackage{amsfonts}
\usepackage{mathtools}
\usepackage[version=4]{mhchem}
\usepackage{siunitx}
\usepackage{longtable,tabularx}
\usepackage{gensymb}
\usepackage{subfigure}
\usepackage{subcaption}
\usepackage{comment}
\usepackage{cite}
\usepackage{authblk}
\usepackage{hyperref}

\newcommand{\origri}[1]{\pmb{\tilde{r}}_{#1}}
\newcommand{\origvi}[1]{\pmb{\tilde{v}}_{#1}}
\newcommand{\origvperi}[1]{\pmb{\tilde{v}}_{{#1}_\perp}}

\newcommand{\Si}[1]{\tilde{S}_{#1}}
\newcommand{\mi}[1]{\tilde{m}_{#1}}

\newcommand{\Dr}{\Delta r}
\newcommand{\mumd}{(\mu-\Dr)}
\newcommand{\thetamd}{\Theta(\mu-\Dr)}
\newcommand{\thetamdtag}{\Theta'(\mu-\Dr)}
\newcommand{\lri}[1]{\pmb{\lambda}_{r_{#1}}}
\newcommand{\lvi}[1]{\pmb{\lambda}_{v_{#1}}}
\newcommand{\lvperi}[1]{\pmb{\lambda}_{{v_{#1}}_\perp}}
\newcommand{\lmu}{\lambda_{\mu}}
\newcommand{\lridot}[1]{\dot{\pmb{\lambda}}_{r_{#1}}}
\newcommand{\lvidot}[1]{\dot{\pmb{\lambda}}_{v_{#1}}}

\newcommand{\statevec}{x}

\newcommand{\origCl}[1]{\tilde{C}l_{\alpha_{#1}}}
\newcommand{\origCdz}[1]{\tilde{C}d_{0_{#1}}}
\newcommand{\origCd}[1]{\tilde{C}d_{\alpha^2_{#1}}}
\newcommand{\origui}[1]{\tilde{u}_{#1}}
\newcommand{\density}{\tilde{\rho}}
\newcommand{\umaxi}[1]{\tilde{U}_{#1}}

\newcommand{\Rturni}[1]{\tilde{R}_{#1}}
\newcommand{\Tturni}[1]{\tilde{T}_{#1}}
\newcommand{\Vinit}[1]{\tilde{V}_{0{#1}}}

\newcommand{\Cdz}[1]{Cd_{0_{#1}}}
\newcommand{\Cd}[1]{Cd_{\alpha^2_{#1}}}

\newcommand{\ri}[1]{\pmb{r}_{#1}}
\newcommand{\vi}[1]{\pmb{v}_{#1}}
\newcommand{\ridot}[1]{\dot{\pmb{r}}_{#1}}
\newcommand{\vidot}[1]{\dot{\pmb{v}}_{#1}}

\newcommand{\vperi}[1]{\pmb{v}_{{#1}_\perp}}
\newcommand{\ui}[1]{u_{#1}}

\newcommand{\secref}[1]{Sec.~\ref{#1}}

\newcommand{\figref}[1]{Fig.~\ref{#1}}

\begin{document}

\title{Flyby Distance Pursuit for Guarding a Target \\
with an Inferior Guard}

\author{Navot Israeli \footnote{Researcher, WALES Ltd., \href{mailto:navot.israeli@gmail.com}{navot.israeli@gmail.com},(Corresponding author).}}
\affil{WALES Ltd., Ramat Gan, Israel, 5252226}
\author{Vitaly Shaferman \footnote{Associate Professor, The Stephen B. Klein Faculty of Aerospace Engineering, Technion, \href{mailto:vitalysh@technion.ac.il}{vitalysh@technion.ac.il}.} and Oded Golan\footnote{Full Research Fellow, The Stephen B. Klein Faculty of Aerospace Engineering, Technion, \href{mailto:oded.golan@technion.ac.il}{oded.golan@technion.ac.il}.}}
\affil{Technion, Israel Institute of Technology, Haifa, Israel, 3200003}

%\institute{Navot Israeli, Corresponding author \at
%             WALES Ltd.\\
%             Ramat Gan, Israel, 5252226 \\
%             navot.israeli@gmail.com
%           \and
%              Vitaly Shaferman  \at
%              Technion, Israel Institute of Technology \\
%              Haifa, Israel, 3200003 \\
%              vitalysh@technion.ac.il
%            \and
%              Oded M. Golan  \at
%              Technion, Israel Institute of Technology \\
%              Haifa, Israel, 3200003 \\
%              oded.golan@technion.ac.il  
%}

%\date{Received: date / Accepted: date}
%The correct dates will be entered by the editor.

\maketitle

\begin{abstract}
Guarding a target against a fast Attacker with a slower Guard is posed as a differential game. Both players follow nonlinear, thrust-free, aerodynamic motion in the plane. The capturability problem inherent in this scenario is addressed by formulating the game such that the Attacker is constrained to terminate at its Target, and the game value is the flyby distance. The latter is made accessible at termination by introducing an auxiliary state variable that records the players' minimum separation distance throughout their flight. This novel combination of termination condition and auxiliary state removes a structural limitation of current pursuit games and opens a new class of solvable problems. It also enables us to embed attacker mission constraints directly. Taking an indirect approach, we formulate the game's Two-Point Boundary Value Problem, solve it numerically, and identify several solution types. When the Attacker's only constraint is to reach the Target, he can do so without being captured. These types of solutions, however, are less practical due to long flight times and slow terminal speeds. Imposing an Attacker terminal speed constraint yields simpler solutions and reveals the trade-off between terminal speed and flyby distance. We demonstrate that in such a setting, capture is possible in a relevant region of the parameter space.  
\end{abstract}

%\keywords{Differential Games \and Guarding a Target \and Capturability \and Slow Guard \and Aerodynamics}
%\subclass{49N75  \and 49N70}
%49N75=Pursuit and evasion games, 49N70=Differential games and control

\section{Introduction}
\label{sec:Introduction}
The classic ``Guarding a Target'' (GaT) game \cite{isaacs_differential_1999,pachter_differential_2017} consists of an Attacker (Evader) that seeks to approach a Target, and a Guard (Pursuer) that seeks to catch it as far from the Target as possible. The original problem was defined with simple motion, but players' dynamics can, in principle, be made more realistic. It is a zero-sum differential game (often termed a pursuit or pursuit-evasion game in the defense context) that assumes both players have conflicting objectives, play optimally, and have full information at all times. Game solutions are optimal Attacker and Guard trajectories that represent an equilibrium in the sense that if either side deviates, it does not gain. Differential games in general, and GaT in particular, are therefore a natural framework for the study of pursuit-evasion strategies in defense scenarios. 
Despite the elegance and promise of this approach, applications in realistic air and missile-defense scenarios have been limited.  

The use of differential games was pioneered by Isaacs \cite{isaacs_differential_1999} and has since become a well-established theory \cite{bryson_applied_1975,basar_dynamic_1999}. A differential game is defined by the players' dynamics, constraints, and objectives. Via the necessary optimality conditions, this leads to a Two-Point Boundary Value Problem (TPBVP), a differential equation for the problem's state and co-state variables, with boundary conditions at the initial and final times. 
One reason for the limited application of differential games in realistic scenarios is that the TPBVP is often difficult to solve. Analytical solutions exist only for simple toy problems and linear systems. Numerical solutions are also quite challenging. For this reason, the application of differential games to air and missile defense was largely limited to problems that could be linearized. Endgame problems often justify linearization, and examples of differential games applications to endgame problems can be found in \cite{ben-asher_advances_1998,gutman_applied_2005,balakrishnan_advances_2013}. 

A second limitation, specific to GaT (and other problems with no predefined termination time), is that it is applicable only in the Usable Part of the terminal capture surface \cite{basar_dynamic_1999}. Namely, only in scenarios where capture is guaranteed by Guard superiority. The limitation arises by construction because Capture is chosen as the game-termination condition. If the Attacker is capable of escape, the termination condition is never satisfied, and the game admits no solutions. 

GaT is therefore inapplicable in scenarios where capture is not guaranteed by the players' intrinsic capabilities. But the fact that GaT is inapplicable does not necessarily mean an Attacker's success. For example, a fast Attacker may lose its velocity advantage on the ingress due to drag forces and be captured by an intrinsically inferior Guard. Note that by giving this simple example, we have introduced a subtle point. When we refer to ``ingress'', we implicitly assume that the Attacker will eventually reach the target. He cannot simply disengage and escape. This fact is absent from GaT, which is agnostic to whether the Attacker is captured on a Target ingress trajectory (In Isaacs’ original simple-motion formulation, this distinction is immaterial, whereas under realistic aerodynamic flight dynamics it may fundamentally alter the engagement outcome). In other words, in our simple example, the Attacker is given a modified mission: reach the Target without being captured. Consequently, the Guard's mission is modified to catch the Attacker on the Target ingress trajectories. These modified missions may be viable to either side, and their success or failure is determined not solely by intrinsic capabilities.

The main contribution of our work is the formulation of a new pursuit game that implements these modified Guard and Attacker missions. We termed it Flyby Distance Pursuit (FBDP). As the name suggests, the game value in FBDP is the players' flyby distance, i.e., the minimum separation throughout their flight. This value is made accessible at game termination by introducing an auxiliary state variable that records the minimum separation. FBDP terminates when the Attacker arrives at the Target, which is the only well-defined event (excluding game start) at our disposal\footnote{One may suggest using flyby as a termination point, but flyby is not a well-defined event because it cannot be recognized in real time. Players' separation may evolve through several minima (see examples in \ref{sec:TypeASolutions}), and the global minimum can only be identified once the Target is reached.}. This terminal condition is convenient because it allows us to introduce terminal constraints that are key to the Attacker's mission (e.g., reaching the Target or maintaining sufficient speed) and limit its evasion abilities. The game admits solutions with arbitrary flyby distances, which can be interpreted as either Guard or Attacker success when compared to an appropriate capture distance. 

Our second contribution is game solutions. The TPBVP arising from FBDP is nonlinear and will probably defy analytic solution attempts. Numerical solutions, however, are possible through the use of CNLP methods that combine Collocation and Non Linear Programming (NLP) solvers \cite{betts_practical_2010}. 
Collocation enables efficient time discretization by transcribing the differential equations into algebraic constraints relating the state variables and controls on a sparse time grid. The problem is then passed to an NLP solver that searches for a satisfying variable assignment.

Horie and Conway \cite{horie_optimal_2006} and Carr et al. \cite{carr_solution_2018} used a CNLP optimal control solver to solve pursuit games. They developed hybrid semi-direct CNLP schemes that perform a direct one-player optimization over the other player's indirect solution. Other authors have suggested iterative schemes in which direct optimization is applied alternately to one player while the other is held fixed \cite{raivio_capture_2001}.  ``Direct'' in this context means optimization via direct manipulation of the problem's controls, while ``indirect'' means optimization by satisfying the necessary conditions (i.e., the TPBVP solution). The appeal of direct optimization is that it eliminates the need for co-state variables, thereby reducing the problem dimension by half. It also eliminates the need to provide an initial guess for the co-state solution, which is often non-intuitive. Direct methods, however, are only appropriate for one-sided optimization, where the objective function is nonconflicting. The hybrid methods developed by these authors were meant to bypass this difficulty, at the cost of resorting to an iterative solver. Following the above examples, we adopt the CNLP approach but insist on a purely indirect solution. We use the optimal control solver PSOPT \cite{becerra_solving_2010,PSOPTman} to solve FBDP TPBVP as is.
\footnote{Although PSOPT is an optimal control solver (designed for optimizing a cost functional of a controlled dynamic system), we only use it as a TPBVP solver with no optimization.}

The rest of the paper is organized as follows. The FBDP is defined in \secref{sec:FlybyDistancePursuit}. The game is formulated for two players who follow thrust-free, aerodynamic motion in a constant altitude plane. Working in a 2D plane is the main simplification that we make. We justify it by the fact that it contains the problem's main ingredients (large-scale maneuvers with energy loss due to drag) and by the insight it provides for solving a 3D problem in the future.  
The FBDP is numerically solved, and the solutions are presented in \secref{sec:FBDPSolutions}. We focus on a head-on engagement of a fast Attacker by a slower Guard. Two types of Attacker terminal constraints are considered. 

In \secref{sec:CATSSolutions} we present solutions for cases where the Attacker is conscious of his terminal speed in order to maintain his kinematic advantage throughout the game. We do this by constraining the Attacker to reach the Target with specific terminal speeds. Varying the terminal speed constraint, we thus draw the Attacker's tradeoff between flyby distance and terminal speed. Our results show that reasonable terminal speed constraints lead to capture. Namely, defense is possible despite the Attacker's kinematic advantage.

\secref{sec:FATSSolutions} is devoted to solutions with free Attacker's terminal speed. These solutions are generally impractical due to long flight times and low terminal speeds. 
They are, however, valuable because they expose multiple local pursuit–evasion solution structures and illustrate the complexity of the underlying game.
We identified four distinct solution types with free Attacker's terminal speeds. In contrast to the constrained terminal speed case, the Attacker can evade the Guard and reach the Target with large flyby distances.    
Some of the free Attacker's terminal speed solutions that we found may correspond to solutions to a degenerate version of the problem, in which the Attacker is allowed only maximal or zero control. This degenerate problem is studied in Appendix \ref{sec:AppendixBZB}. We find it useful for understanding the solution structure of the FBDP. Appendix \ref{sec:AppendixBZB} may also be useful for supplying initial guesses for the CNLP solver. 

Our conclusions are given in \secref{sec:Conclusions}.

\section{Flyby Distance Pursuit}
\label{sec:FlybyDistancePursuit}
We begin by formulating the problem as a zero-sum differential game between Attacker and Guard. Both players move in a constant altitude $(\tilde{x},\tilde{y})$ plane, 
with air density $\density$. Their positions and velocities are denoted, respectively, by 
\begin{equation}
\label{eq:player_states}
\origri{i} =\left[\tilde{x}_i, \tilde{y}_i\right]^T, \quad 
\origvi{i}=\left[\tilde{v}_{x_i}, \tilde{v}_{y_i}\right]^T
\end{equation}
where $i \in \{a,g\}$  stand for Attacker/Guard respectively. We decorate these variables and parameters with $\tilde{\Box}$ because we will soon re-scale the problem and prefer the re-scaled expressions to be concise.  

The players obey equations of motion of a point mass, thrust-free, instantaneous Angle of Attack (AoA) control:  
\begin{equation}
\label{eq:eom}
\begin{split}
\frac{d\origri{i}}{d\tilde{t}}=&\origvi{i} \\
\frac{d\origvi{i}}{d\tilde{t}} =&\frac{\Si{i}\cdot \density \cdot ||\origvi{i}||}{2\cdot \mi{i} } \cdot \bigg[\origCl{i}\cdot \origui{i}\cdot \origvperi{i}-(\origCdz{i}+\origCd{i}\cdot \origui{i}^2)\cdot \origvi{i}\bigg]
\end{split}
\end{equation}

The initial conditions are 
\begin{equation}
\label{eq:initial_condition}
\origri{i}(\tilde{t}=0)=\origri{i_0}, \ \ \origvi{i}(\tilde{t}=0)=\origvi{i_0} 
\end{equation}
and the final conditions will be specified later as part of the formulation of the game. 
$\Si{i}$, $\mi{i}$, $\origCl{i}$, $\origCdz{i}$, and $\origCd{i}$ are the players' reference area, mass, lift slope, parasitic drag, and induced drag coefficients.   $\origui{i}$ are the players' AoA controls in the horizontal plane. They are limited by $|\origui{i}|\leq \umaxi{i}$. $\origvperi{i}=\left[-\tilde{v}_{y_i},  \tilde{v}_{x_i}\right]^T$ are the perpendicular velocity vectors. In the sequel we will use $\pmb{z}_\perp\equiv\left[-z_y,z_x\right]^T$ for any two dimensional vector $\pmb{z}=\left[z_x, z_y\right]^T$.
Note that the players' radii of maximal turn (i.e., minimal radii) $\Rturni{i}=\frac{2\cdot \mi{i}}{\Si{i}\cdot \density \cdot \origCl{i} \cdot \umaxi{i}}$  are time independent. 
Adopting $\Rturni{a}$ and $\Tturni{a}=\frac{\Rturni{a}}{\Vinit{a}}$ as our length and time units, we move to the scaled position   $\ri{i}=\frac{\origri{i}}{\bar{R}_a}$, velocity $\vi{i}=\frac{\origvi{i}}{\Vinit{a}}$ and time  $t=\frac{\tilde{t}}{\Tturni{a}}$ variables. $\Vinit{a}=||\origvi{a}(0)||$ here is the Attacker's initial speed. The equations of motion Eq. \ref{eq:eom} are rewritten as 
\begin{equation}
\label{eq:scaled_eom}
\ridot{i}=\vi{i},\;
\vidot{i}=\zeta_i\cdot 
||\vi{i}|| \cdot \bigg[ \ui{i}\cdot\vperi{i}-(\Cdz{i}+\Cd{i}\cdot \ui{i}^2)\cdot \vi{i}\bigg]
\end{equation}
where $\dot{\Box}\equiv\frac{d\Box}{dt}$ and
\begin{equation}
\label{eq:scaled_par}
\Cdz{i}=\frac{\origCdz{i}}{\origCl{i}\cdot \umaxi{i}}, \quad \Cd{i}=\frac{\origCd{i} \cdot \umaxi{i}}{\origCl{i}}, \quad
\ui{i}=\frac{\origui{i}}{\umaxi{i}}, \quad \zeta_a\equiv1, \quad \zeta_g=\frac{\Rturni{a}}{\Rturni{g}}
\end{equation}
 
The scaled controls $\ui{i}$ are bounded by $|\ui{i}|\leq1$. We are thus left with only 5 parameters (compared to the initial 13).

Using these scaled variables, the game state vector is given by 
\begin{equation}
\label{eq:game_state}
\pmb{\statevec} = \left[  \ri{a}^T, \ri{g}^T ,\vi{a}^T, \vi{g}^T, \mu\right]^T
\end{equation}
where $\mu$ is an auxiliary variable that records the players' minimal separation.

\subsection{Flyby Distance Variable}
\label{sec:FlybyDistanceVariable}
The flyby distance variable obeys the following dynamics:
\begin{equation}
\label{eq:mudynamics}
\dot{\mu}=-\frac{\thetamd \cdot (\mu-\Dr)}{\tau}, \quad 
\Dr=\sqrt{||\ri{a}-\ri{g}||^2+\delta^2}
\end{equation}
where $\Dr$ is the players' distance with the cusp smoothed on a small scale $\delta$, and $\tau$ is a small time constant. $\Theta\left(\cdot\right)$ is a step function which is $\delta$ smoothed as well. We use the piecewise continuous expression 

\begin{equation}
\Theta(x)=\left \{\begin{array}{cc}
     1,& \delta<x  \\
     1-\frac{(x-\delta)^2}{2\cdot \delta^2},& 0<x\leq \delta \\
     \frac{(x+\delta)^2}{2\cdot \delta^2},& -\delta<x\leq 0 \\
     0,& x<-\delta
\end{array}
\right. 
\end{equation}
which has a continuous first derivative. 

By this construction, the flyby distance variable $\mu$ converges exponentially to $\Delta r$ when the players' separation decreases below its former minimal value and remains constant otherwise. The scales $\tau$, $\delta$ serve as small parameters that define the temporal and spatial resolution of the process. Smoothing by $\delta$ here simply means that we are not interested in flyby distances smaller than $O(\delta)$, which may loosely be interpreted as the capture radius (in our results, see \figref{fig:CATS0.4_mu}, %Oded 14/4 this means that the capture must be less than 2*delta
we find that $\mu(t_f)\approx2\cdot\delta$ means almost zero actual flyby distance). With a constant closing speed $\dot{\Delta r}$, $\mu$ converges exponentially to a steady state that lags behind $\Dr$ a distance $\tau\cdot \dot{\Delta r}$. The time constant $\tau$ should therefore be on the scale of $\frac{\delta}{\dot{\Delta r}}$.    

\subsection{Objective Function}
\label{sec:ObjectiveFunction}
The Attacker and Guard maximize and minimize, respectively, the following Mayer cost:
\begin{equation}
    \label{eq:terminal_cost}
    \mathcal{J}=\mu(t_f)+\phi_{v_a}\cdot ||\vi{a}(t_f)||-\phi_{v_g}\cdot||\vi{g}(t_f)||
\end{equation}
The primary cost is the terminal flyby distance $\mu(t_f)$. To this, we add terms that favor speed conservation. These terms are necessary because the terminal flyby distance is determined at some intermediate time $t_{fb} \leq t_f$. Valuing the flyby distance alone provides no guidance (especially for the Guard) for $t_{fb}\leq t \leq t_f$. 
Adding small $\phi_{v_a},\phi_{v_g}$ regularization terms removes this degeneracy without significantly altering the game.
Adding significant $\phi_{v_a}, \phi_{v_g}$ costs may be useful for studying the players' trade-off between flyby distance and final speeds. This, however, is more conveniently achieved by setting concrete terminal speed constraints (see \ref{sec:BoundaryConditions}). In the sequel, $\phi_{v_a},\phi_{v_g}$ will only be used as small regulatory terms.

\subsection{Hamiltonian and Co-states Dynamics}
\label{sec:Hamiltonian}
Following \cite{basar_dynamic_1999} we write the game Hamiltonian
\begin{equation}
\label{eq:Hmiltonial}
\mathcal{H}=\lri{a}^T \ridot{a} +\lri{g}^T \ridot{g} +\lvi{a}^T \vidot{a} +\lvi{g}^T \vidot{g} + \lmu \dot{\mu}
\end{equation}
which defines the co-states 
\begin{equation}
    \label{eq:Costatesdef}
    \pmb{\lambda} = \left[\lri{a}^T , \lri{g} ^T, \lvi{a}^T, \lvi{g}^T, \lmu\right]^T
\end{equation}

Co-states dynamics obey $\dot{\pmb{\lambda}}=-\left( \frac{\partial \mathcal{H}}{\partial \pmb{x} }\right)$, which expands to
\begin{equation}
\label{eq:costate_dynamics}
\begin{split}
\lridot{i}=&-\frac{\left(\ri{i}-\ri{\bar{i}}\right)}{\Dr}\cdot \frac{\lmu}{\tau}\cdot \left[\thetamdtag\cdot \mumd + \thetamd\right] \\
\lvidot{i}=&-\lri{i} 
+\zeta_i\cdot
\left\{ 
\left[
\left(\Cdz{i}+\Cd{i}\cdot \ui{i}^2 \right)
\cdot 
\left(\lvi{i}^T\vi{i}\right)-\ui{i}\cdot \left(\lvi{i}^T\vperi{i}\right)\right]\cdot\frac{\vi{i}}{||\vi{i}||} \right. \\
&+\left. ||\vi{i}||\cdot\left[ \ui{i}\cdot\lvperi{i} +   \left(
\Cdz{i}+\Cd{i}\cdot \ui{i}^2 
\right)  \cdot \lvi{i}\right] \right \} \\
\dot{\lambda}_{\mu}=&\frac{\lmu}{\tau}\cdot 
\left[
\thetamdtag\cdot\mumd + \thetamd
\right] 
\end{split}
\end{equation}
where $\bar{i}=\{a,g\}\setminus i$ is the opponent of $i$.

\subsection{Optimal Controls}
\label{sec:OptimalControls}
The controls $\ui{a}$ and $\ui{g}$ enter the Hamiltonian in Eq. \ref{eq:Hmiltonial} only through the terms $\lvi{a}^T\vidot{a}$ and $\lvi{g}^T\vidot{g}$ respectively. Isaacs condition (see \cite{basar_dynamic_1999} p. 428) that $\underset{|\ui{a}|\leq1}{\mathrm{max}} \left[\underset{|\ui{g}|\leq 1}{\mathrm{min}}(\mathcal{H})\right]=\underset{|\ui{g}|\leq 1}{\mathrm{min}} \left[\underset{|\ui{a}|\leq 1}{\mathrm{max}}(\mathcal{H})\right]$ therefore holds and the optimal controls are given by

\begin{equation}
\label{eq:optimal_controls}
\begin{split}
\ui{a}^*=\underset{|\ui{a}|\leq 1}{\mathrm{argmax}}\left(\lvi{a}^T\vidot{a}\right) \\
\ui{g}^*=\underset{|\ui{g}|\leq 1}{\mathrm{argmin}}\left(\lvi{g}^T\vidot{g}\right) \\
\end{split}
\end{equation}

Therefore, the Attacker's acceleration is attracted to his velocity co-state direction, and the Guard's acceleration is repelled from his. Using Eq. \ref{eq:scaled_eom} and omitting controls independent terms, we get

\begin{equation}
\ui{a}^*=\underset{|\ui{a}|\leq 1}{\mathrm{argmax}}\left(\ui{a}\cdot\lvi{a}^T\vperi{a}-\Cd{a}\cdot \ui{a}^2\cdot \lvi{a}^T\vi{a}\right) 
\end{equation}
which depends on  $\vi{a}$ and its co-state  $\lvi{a}$ only through the angle 
\begin{equation}
\varphi_a=atan2\left( \lvi{a}^T \vperi{a},\lvi{a}^T \vi{a}\right)    
\end{equation}
between them. Replacing 
\begin{equation}
\lvi{a}^T\vperi{a}=||\vi{a}||\cdot ||\lvi{a}||\cdot \sin{\varphi_a}, \quad  
\lvi{a}^T \vi{a}=||\vi{a}||\cdot ||\lvi{a}||\cdot \cos{\varphi_a}
\end{equation}
and canceling the $||\vi{a}||\cdot ||\lvi{a}||$ factors, we are left with

\begin{equation}
\label{eq:optimal_controls2}
\ui{i}^*=\underset{|\ui{i}|\leq 1}{\mathrm{argmax}}\left[\sin{\varphi_i} \cdot \ui{i} - \Cd{i}\cdot\cos{\varphi_i}\cdot \ui{i}^2\right]\\
\end{equation}
which is applicable to the Guard as well, provided that we define
\begin{equation}
\varphi_g=atan2\left( -\lvi{g}^T \vperi{g},-\lvi{g}^T \vi{g}\right)
\end{equation} 
as the angle between $-\lvi{g}$ and $\vi{g}$.
Equation \ref{eq:optimal_controls2} is solved by
\begin{equation}
\label{eq:optimal_controls3}
\ui{i}^*=\left \{\begin{array}{cc}
    1,& \varphi_i>\bar{\varphi}_i \\
     \frac{\tan{\varphi_i} }{2 \cdot \Cd{i}} ,& |\varphi_i|\leq \bar{\varphi}_i  \\
     -1,& \varphi_i<-\bar{\varphi}_i
\end{array}
\right. 
\end{equation}
with $\bar{\varphi}_i=\tan^{-1} \left(2\cdot \Cd{i} \right)$.
Note that the induced drag coefficient $\Cd{i}$ introduces a region $\left(-\bar{\varphi}_i,\bar{\varphi}_i\right)$ of moderated control. It has a stabilizing effect on the equations of motion Eq. \ref{eq:scaled_eom} integration, which may otherwise be numerically unstable. This is especially important at the termination time, where the boundary condition Eq. \ref{eq:transversality_bc} dictates that the velocities will be aligned with their co-states. 

\subsection{Boundary Conditions}
\label{sec:BoundaryConditions}
A solution of the game differential equations defined in Eq.  \ref{eq:scaled_eom} and Eq. \ref{eq:costate_dynamics} requires a total of 19 boundary conditions (8 kinematic state variables, 8 kinematic co-state variables, 1 flyby distance, 1 flyby distance co-state, 1 free termination time). 

Eight boundary conditions are determined by the initial state of Eq. \ref{eq:initial_condition}:
\begin{equation}
\label{eq:kinematic_bc}
\begin{split}
\ri{a}(t=0)=\ri{a_0}, \ \ \vi{a}(t=0)=\vi{a_0} \\
\ri{g}(t=0)=\ri{g_0}, \ \ \vi{g}(t=0)=\vi{g_0} 
\end{split}
\end{equation}

A ninth condition specifies the initial flyby distance value: 
\begin{equation}
\label{eq:flyby_bc}
\mu(t=0)=\Dr(t=0)-\tau\cdot \dot{\Delta r}(t=0)
\end{equation}
This boundary condition puts $\mu$ in a steady state distance from $\Dr$ when the two players are initially moving closer with a constant closing speed (i.e., $\dot{\Delta r} <0$). As $\mu(t)$ converges to $\Dr(t)-\tau\cdot \dot{\Delta r}(t)$ exponentially, this boundary condition has only a transient effect and does not change the overall solution significantly.

As mentioned before, the game terminates when the Attacker reaches the Target. Namely, 
\begin{equation}
\label{eq:termination_condition}
\ri{a}\left( t_f \right)=\ri{target}
\end{equation}
We may also wish to apply an Attacker Terminal Speed (ATS) constraint \footnote{An inequality constraint here would have been more natural, but is harder to solve.}. 
\begin{equation}
\label{eq:termination_speed_condition}
||\vi{a}\left( t_f \right)||=V_{a_{terminal}}
\end{equation}

These terminal conditions define a termination manifold upon which the vector
\begin{equation}
\pmb{\Gamma} = \left[\ri{a}^T(t_f)-\ri{target}^T,  ||\vi{a}\left( t_f \right)||-V_{a_{terminal}}\right]^T   
\end{equation}
vanishes.

The transversality conditions are given by the gradient of the terminal cost with respect to the game final state, on the termination manifold:
\begin{equation}
\label{eq:transversality_def}
\pmb{\lambda}(t_f)=\left. \left(\frac{\partial \mathcal{J}}{\partial \pmb{\statevec}(t_f)} \right)^T\right|_{\pmb{\Gamma}=\pmb{0}}
\end{equation}
 This is achieved by replacing $\mathcal{J}$ in Eq. \ref{eq:transversality_def} with the augmented terminal cost 
 \begin{equation}
\label{eq:augmented_cost}
\bar{\mathcal{J}}=\mathcal{J} + \pmb{\kappa}^T\pmb{\Gamma}
\end{equation}
where $\pmb{\kappa}=[\kappa_1,\kappa_2, \kappa_3]^T$ is a vector of arbitrary Lagrange multipliers.  

We get
\begin{equation}
\label{eq:transversality_bc}
\begin{split}
\lri{a}\left(t_f\right)=&\left(\frac{\partial \bar{\mathcal{J}}}{\partial \ri{a}}\right)^T= [\kappa_1 ,\kappa_2]^T  \\
\lri{g}\left(t_f\right)=&\left(\frac{\partial \bar{\mathcal{J}}}{\partial \ri{g}}\right)^T= \bf{0} \\
\lvi{a}\left(t_f\right)=&\left(\frac{\partial \bar{\mathcal{J}}}{\partial \vi{a}}\right)^T=\left(\phi_{v_a}+\kappa_3\right)\cdot \frac{\vi{a}(t_f)}{||\vi{a}(t_f)||}\\
\lvi{g}\left(t_f\right)=&\left(\frac{\partial \bar{\mathcal{J}}}{\partial \vi{g}}\right)^T=-\phi_{v_g}\cdot \frac{\vi{g}(t_f)}{||\vi{g}(t_f)||} \\
\lmu\left(t_f\right)=&\frac{\partial \bar{\mathcal{J}}}{\partial \mu}=1
\end{split}
\end{equation}
where $\bf{0}$ is a vector of zeros with appropriate dimensions. Note that since $\pmb{\kappa}$ is arbitrary, $\lri{a}\left(t_f\right)$ is completely free. In addition, when applying the ATS constraint, the weight $\phi_{v_a}$ is swallowed by $\kappa_3$. $\phi_{v_a}$, therefore, is meaningful only without the ATS constraint. The boundary condition for the Attacker's velocity co-state can thus be summarized by 
\begin{equation}
\label{eq:attacker_lambdav_bc_wfvconst}
\lvi{a}\left(t_f\right)\parallel\vi{a}(t_f) 
\end{equation}
with ATS constraint, and
\begin{equation}
\label{eq:attacker_lambdav_bc_nofvconst}
\lvi{a}\left(t_f\right)=\phi_{v_a}\cdot \frac{\vi{a}(t_f)}{||\vi{a} (t_f)||} 
\end{equation}
without ATS constraint. Note that these two types of boundary conditions are interchangeable.  
Any solution that satisfies the ATS boundary conditions \eqref{eq:termination_speed_condition} and \eqref{eq:attacker_lambdav_bc_wfvconst} also satisfies the free ATS boundary condition \eqref{eq:attacker_lambdav_bc_nofvconst} with $\phi_{v_a}=\frac{\lvi{a}^T\left(t_f\right)\vi{a}\left(t_f\right)}{||\vi{a}\left(t_f\right)||}$

The last boundary condition is \cite{bryson_applied_1975}
\begin{equation}
\label{eq:freetime_bc}
\mathcal{H}\left(t_f\right)=0
\end{equation}

Note that boundary condition $\lmu\left(t_f\right)=1$ together with the EOM in Eq. \ref{eq:costate_dynamics} suggests the approximation
\begin{equation}
\label{eq:approxlmu}
\lmu(t)\approx\left \{\begin{array}{cc}
    e^{\frac{t-t_{fb}}{\tau}}& t\leq t_{fb} \\
     1& t>t_{fb}
\end{array}
\right. 
\end{equation}
Where $t_{fb}$ is the flyby time.

\section{Flyby Distance Pursuit Solutions}
\label{sec:FBDPSolutions}
\subsection{Scenario}
\label{sec:Scenario}
Consider a scenario where the Target and the Guard's launch point are co-located at the origin. The Attacker is detected at a range of  $\approx6$ along the $x$-axis, heading towards the Target. After detection, the Guard needs to make a launch decision. He waits for the Attacker to approach a range $L$ and then launches to a selected azimuth. Prior to the Guard's launch, the Attacker has no reason to turn. Any change of heading on his part will be answered effortlessly by the Guard's launch azimuth. In this state, the Guard launches towards the Attacker, along the positive $x$-axis. The game is thus started at a range $L$, in a head-on geometry. Game solutions from this initial state have a reflection symmetry about the $x$-axis. In the sequel, we focus on solutions in which the Attacker initially plays towards the positive $y$ direction.   

We are interested in a situation in which the Attacker has a speed advantage while the Guard has a smaller (better) turn radius. For simplicity, this is achieved by the following parameter selection. Both players share the same aerodynamic properties, adopted (with some simplifications) from \cite{rizvi_optimal_2015}, with $\origCl{i}=0.75$, $\origCdz{i}=0.012$, $\origCd{i}=0.6$. In reality, they will differ, but we wish to focus the analysis on kinematical aspects. The Guard's smaller turn radius is achieved by giving it a larger maximal AoA of  $\umaxi{g}=30\degree$ while the Attacker is limited by  $\umaxi{a}=15\degree$. In terms of our scaled parameters this is translated to $\Cdz{a}=0.061$, $\Cd{a}=0.209$, $\Cdz{g}=0.031$, $\Cd{g}=0.419$ and $\zeta_g=2$. By our definition of the scaled variables, the initial Attacker's speed is 1. The slower guard has an initial speed of 0.4. Although the Attacker slows down during his flight from the detection point to $L$, he can still produce larger lateral accelerations. For example, with $L=3$, the Attacker's initial maneuver capability is more than twice the Guard's. 

In the general scenario above, we identified several solution types, with and without ATS constraints. Constrained Attacker Terminal Speed (CATS) solutions are more realistic and convey the important result that defense is possible despite the Attacker's kinematic advantage. Solutions without ATS constraints (and small speed constants $\phi_{v_a},\phi_{v_g}$) tend to be impractical due to slow terminal speeds and long flight times. We refer to them as Free Attacker Terminal Speed (FATS) solutions. They are interesting from the pattern formation perspective and may suggest different types of pursuit-evasion strategies. In our workflow, FATS were found first and later used as initial guesses for CATS. Here, we discuss CATS solutions first because of their practical importance.  

All the results presented in the sequel were computed using $\delta,\tau<<1$. The Guard's speed cost weight was set to $\phi_{v_g}=0.37$. An identical Attacker's speed cost weight of $\phi_{v_a}=0.37$ was used for the FATS solutions as well.

\subsection{Constrained Attacker Terminal Speed Solutions}
\label{sec:CATSSolutions}
\begin{figure}[bt!]
\centering
\includegraphics[width=0.7\textwidth]{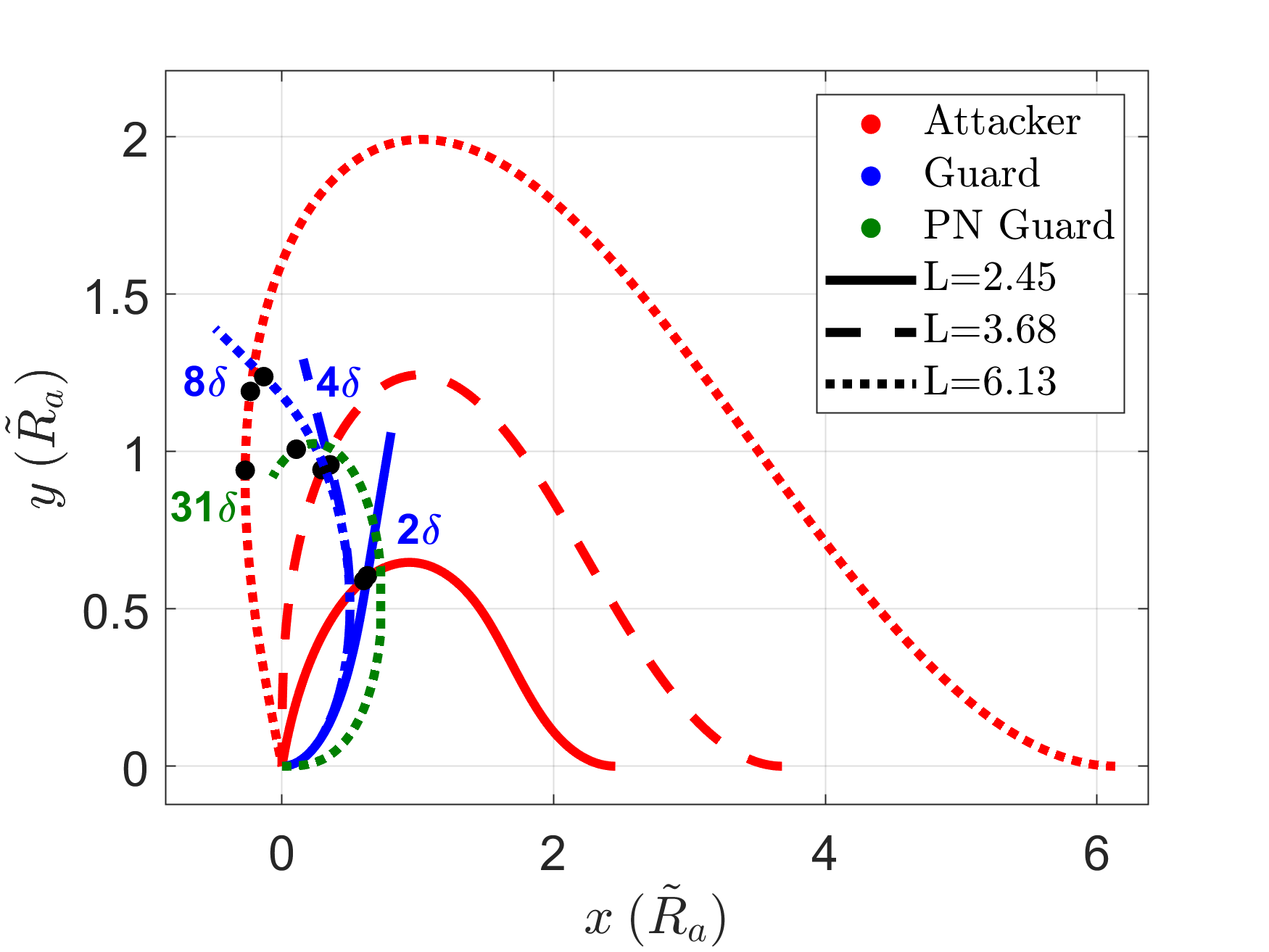}
\caption{Players' trajectories with 0.4 Attacker terminal speed constraint.}
\label{fig:CATS0.4_trajs}
\end{figure}

The Attacker should not arrive at the Target too slow. Approaching the Target at a speed much slower than that of a Guard makes it vulnerable to additional interception attempts. With this in mind, we solved FBDP with different ATS constraint values and from different launch ranges $L$.  
Figure \ref{fig:CATS0.4_trajs} presents players' trajectories from game solutions with an ATS constraint of 0.4, matching the Guard's launch speed. The Attacker is plotted in red while the Guard is in blue. Black markers on the players' trajectories mark the flyby points, and the final value of the flyby distance variable $\mu(t_f)$ is denoted (in $\delta$'s). The Attacker performs an evasion maneuver to the right, flies straight, and then makes a left turn that reaches the Target. The Guard makes a left turn to intercept the Attacker and, once the flyby point is passed, continues straight to conserve his speed. When the Guard launches at a range $L=6.13$, he misses the Attacker by $8\delta$. The flyby distance decreases with launch range, and at $L=2.45$ it shrinks to $2\delta$. 

The Guard's performance in \figref{fig:CATS0.4_trajs} can be appreciated by a comparison with standard guidance laws. To this end, we plotted the green curve in \figref{fig:CATS0.4_trajs} that shows the trajectory of a Proportional Navigation (PN) guided Guard, trying to intercept the $L=6.13$ Attacker trajectory. Compared with the game solution, the PN Guard is slow to apply a maximal turn at first and later does not relax his maximal turn at all (to conserve speed). The end result is a flyby distance of $31\delta$, larger than the $8\delta$ flyby of the game solution. These results were obtained using a navigation constant of $N=4$. A sharper reaction may be elicited by using a larger value of $N$, but $N\approx100$ is required for a good imitation of the early part of the Guard's solution. Even with an impractical $N=100$, the PN Guard performance is inferior to the game solution.              

\begin{figure}[bt!]
\centering    
\subfigure[Players speeds.]{\label{fig:CATS0.4_speed}\includegraphics[width=0.45\textwidth]{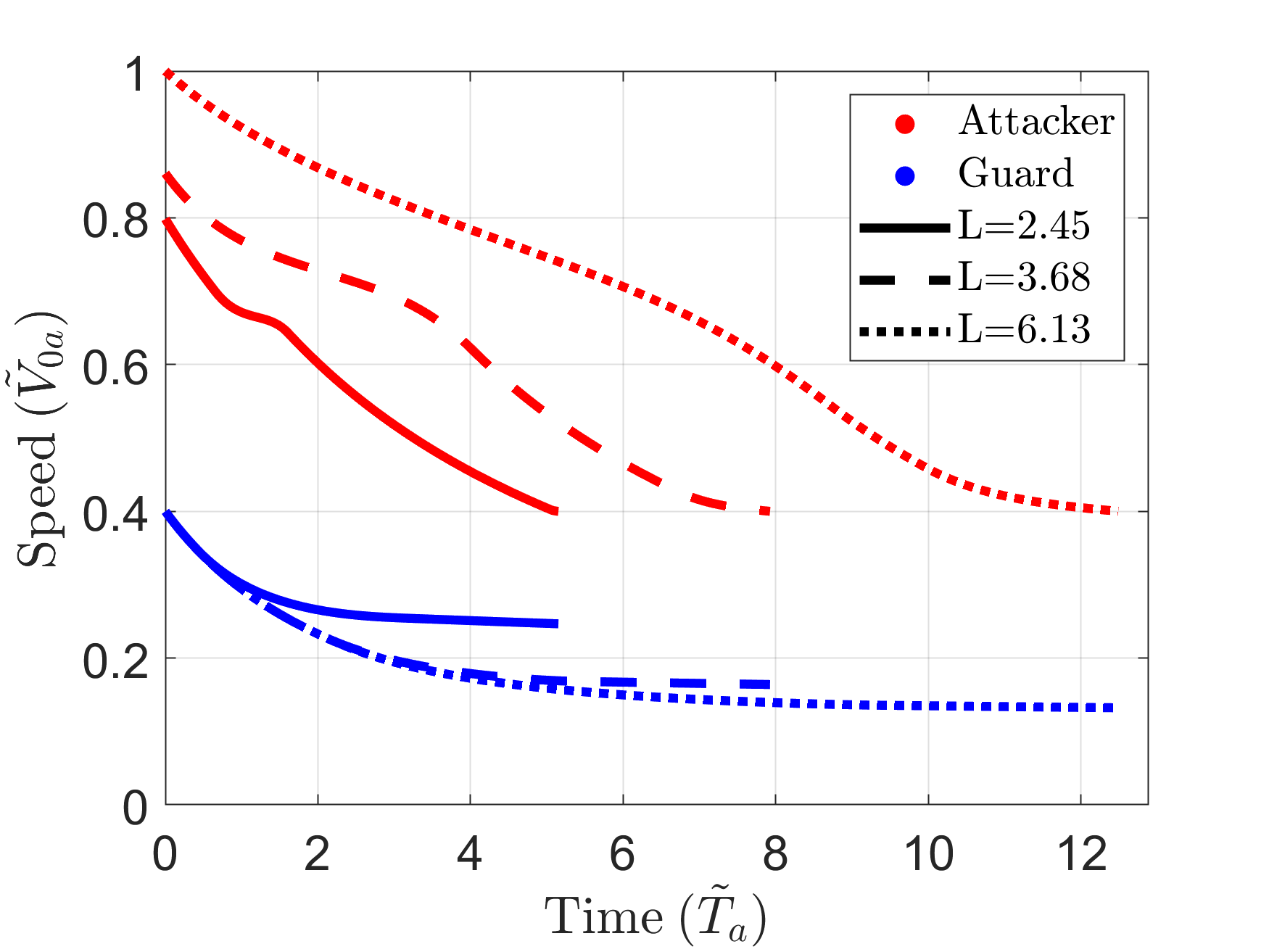}}
\subfigure[Players controls.]{\label{fig:CATS0.4_controls}\includegraphics[width=0.45\textwidth]{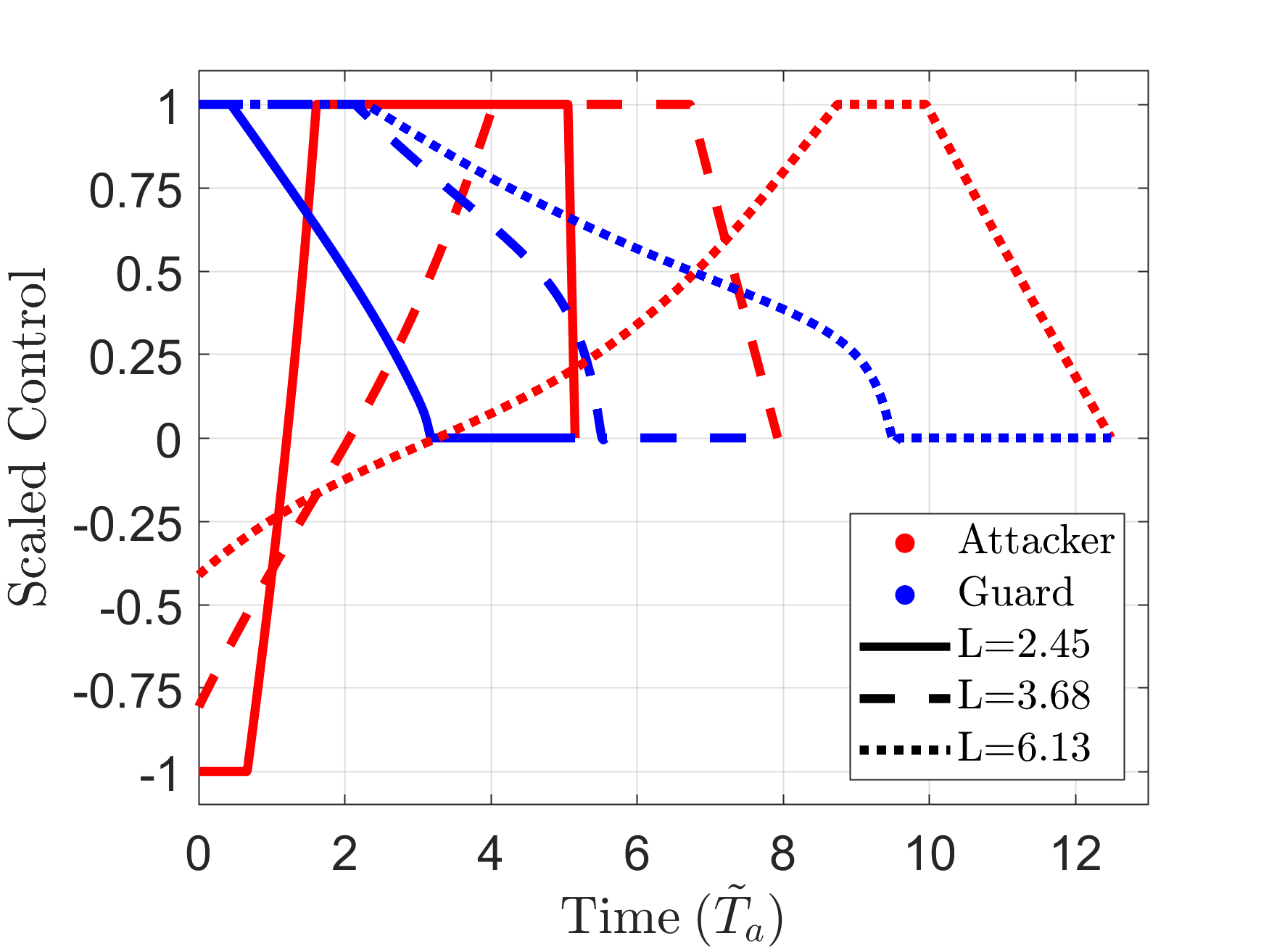}}
\subfigure[Flyby variable.]{\label{fig:CATS0.4_mu}\includegraphics[width=0.45\textwidth]{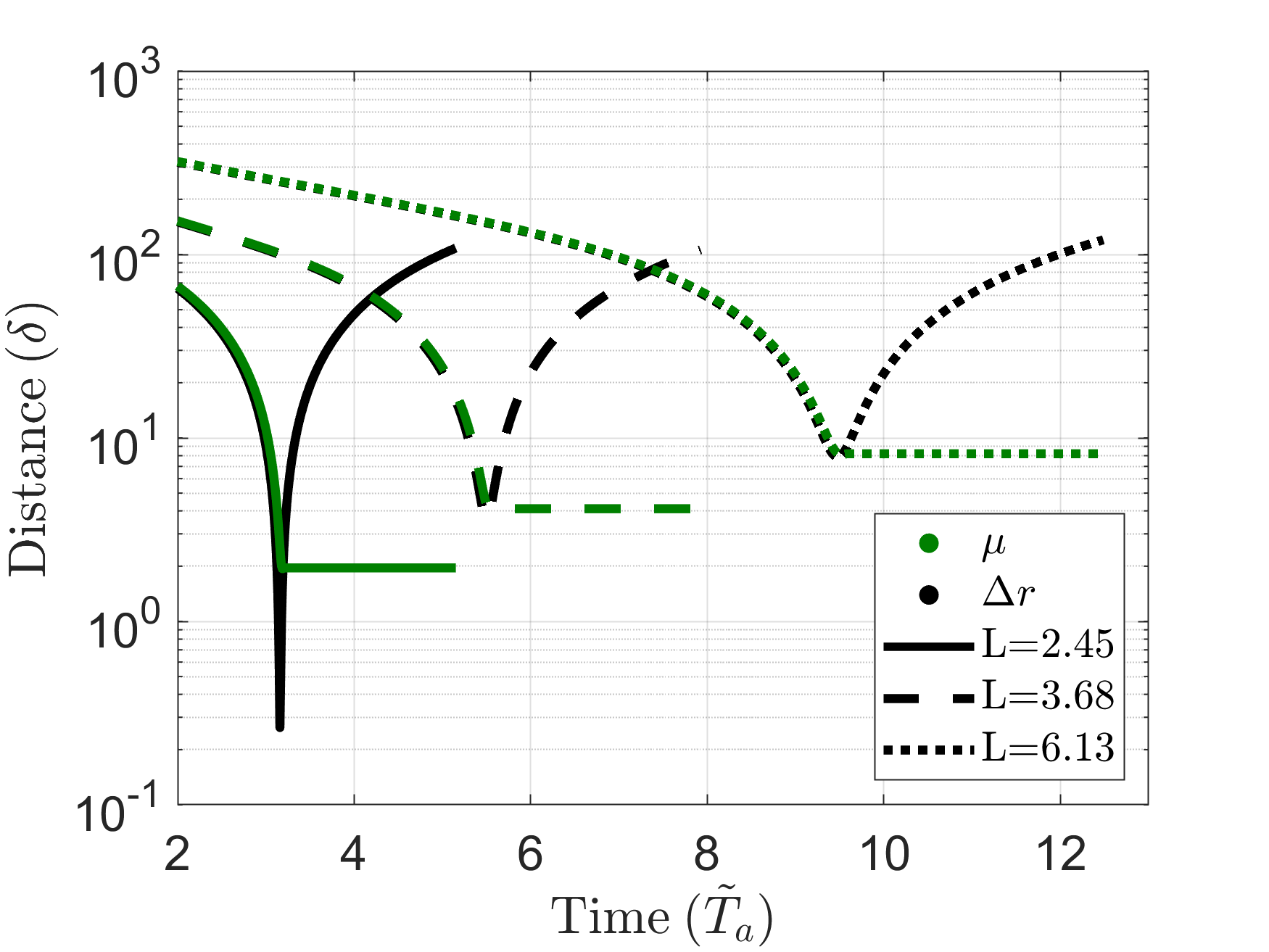}}
\subfigure[Flyby variable co-state.]{\label{fig:CATS0.4_lambdamu}\includegraphics[width=0.45\textwidth]{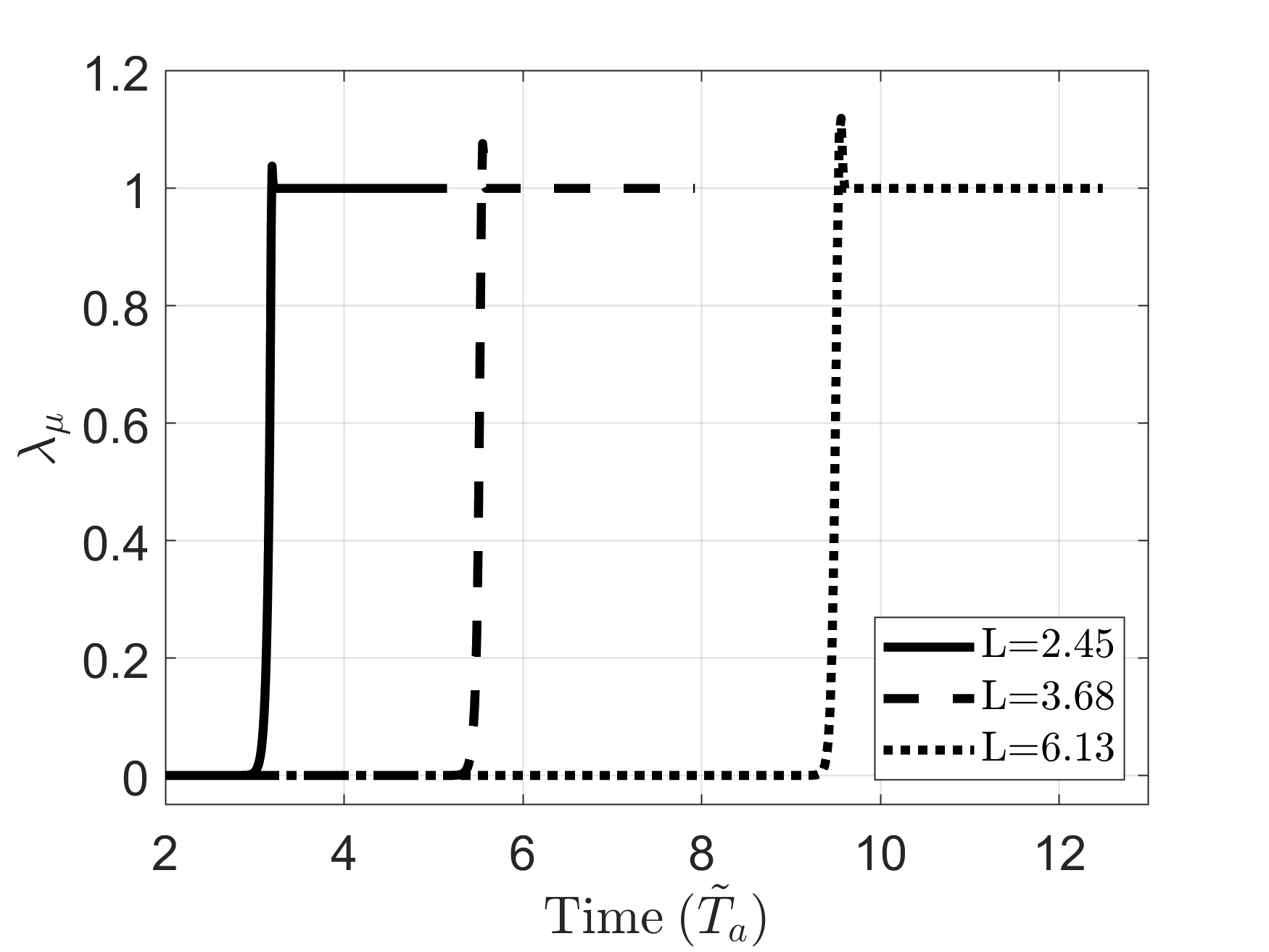}}
\caption{Game solutions with 0.4 Attacker terminal speed constraint.}
\label{fig:CATS0.4_additional_info}
\end{figure}

Figure \ref{fig:CATS0.4_additional_info} presents additional information for the trajectories presented in \figref{fig:CATS0.4_trajs}. Figure \ref{fig:CATS0.4_speed} shows the players' speeds, which decrease considerably with flight time. Figure \ref{fig:CATS0.4_controls} shows the players' scaled controls. At the longest launch range $L=6.13$, the Attacker's control is moderated while the Guard uses a maximal turn at first, reduces it at some point, and then sets it to zero at the flyby time. At shorter launch ranges, the Attacker's control approaches its limiter while the Guard's becomes moderated. At $L=2.45$, the Attacker's control is almost maximal. It maximizes at $L=2.21$ (not shown), and below this range, there are no further solutions because the Attacker cannot shed sufficient speed from a shorter range. Note that both players always terminate with their controls set to zero, in agreement with the optimal control Eq. \ref{eq:optimal_controls3} and the transversality boundary conditions in Eq. \ref{eq:transversality_bc} that align the velocity co-states with the velocity vectors at $t_f$.

Figure \ref{fig:CATS0.4_mu} presents the flyby distance variable $\mu$ (green) and compares it with the actual players' instantaneous separation distance (black). We observe that $\mu$ indeed follows the instantaneous separation until the flyby point and then detaches from it with a value that records the flyby distance. As the solid curves in \figref{fig:CATS0.4_mu} indicate, the recorded value is an upper bound with an accuracy of $\approx 2\delta$. Finally, \figref{fig:CATS0.4_lambdamu} presents the flyby variable costate $\lambda_{\mu}$. As expected from Eq. \ref{eq:approxlmu}, $\lambda_{\mu}$ exhibits a sharp transition from 0 to 1 at the flyby time. The overshoots, however, are not accounted for by the approximation in Eq. \ref{eq:approxlmu} and are introduced by the finite width of the step function $\Theta$. The players' position and velocity co-states are also part of the solution, but we refrain from presenting them, as they are not very informative.  

\begin{figure}[bt!]
\centering
\includegraphics[width=0.7\textwidth]{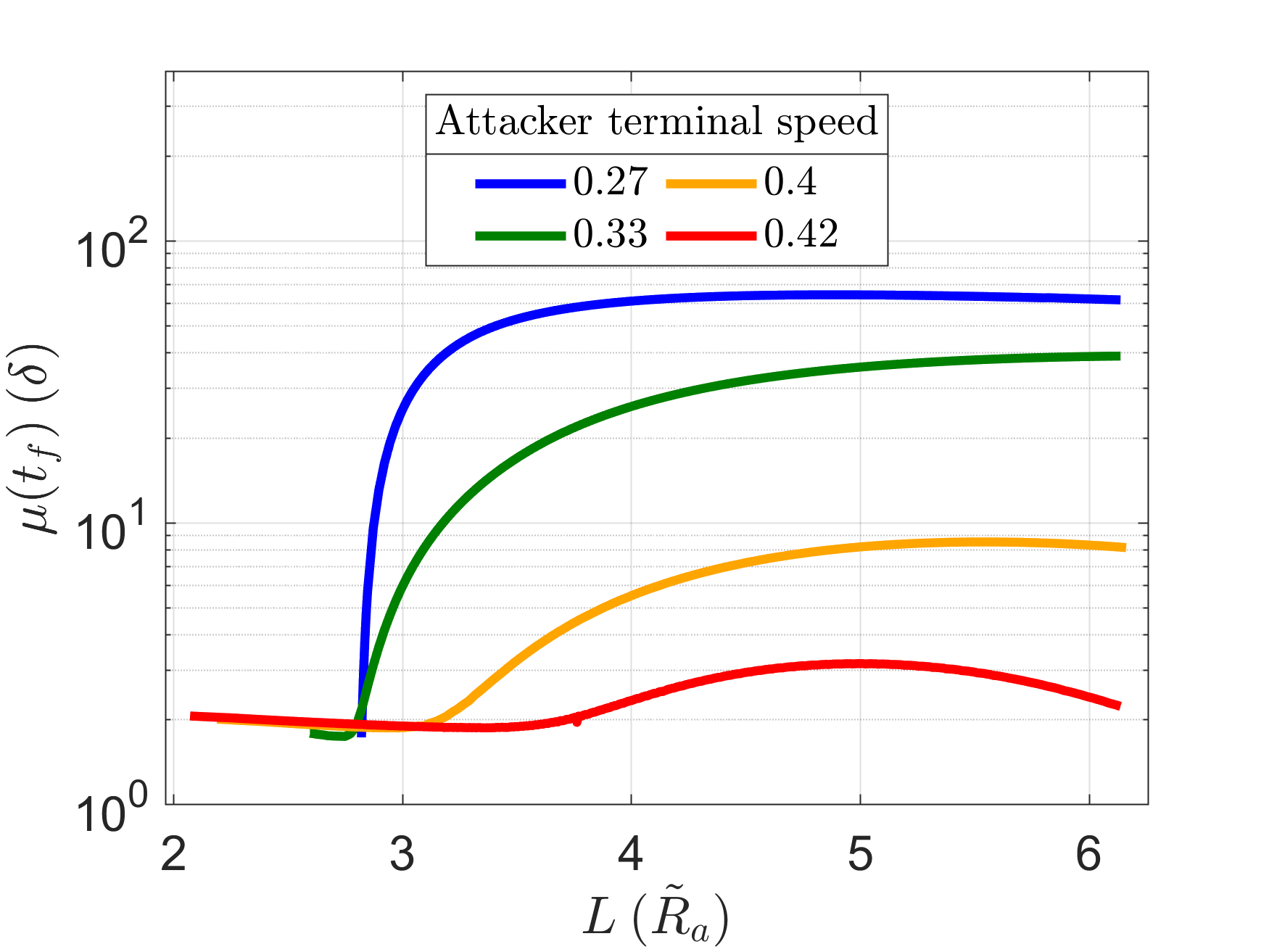}
\caption{Flyby distances with different Attacker terminal speed constraints.}
\label{fig:CATS_flyby_tradeoff}
\end{figure}

Figure \ref{fig:CATS_flyby_tradeoff} presents the Attacker's tradeoff between flyby distance and terminal speed. The different curves show the recorded flyby distance $\mu(t_f)$ as a function of launch range $L$ for different values of the ATS constraint. With ATS constraint of 0.27 (blue), the Attacker can escape with flyby distances on the order of the Attacker's turn radius, unless the Guard launches at a range of $L\approx2.8$, where $\mu(t_f)$ drops below $2\delta$. As we mentioned earlier, $\mu(t_f)$ is an upper bound for the actual flyby distance with an accuracy of $2\delta$. $\mu(t_f)\leq2\delta$, therefore, means that the Attacker is captured. At this point, the blue curve terminates because the Attacker cannot shed enough speed when starting to maneuver at shorter ranges. Flyby distances are reduced with increasing ATS constraint. ATS constraint of 0.33 (green) still results in large flyby distances at large launch ranges and capture at $L\approx2.8$. With ATS constraint of 0.4, the flyby distance is always below $10\delta$ and drops to $2\delta$ at $L\approx3 $. Increasing the ATS constraint to 0.42 (red) extends the capture launch range to 3.7. The implications of these results are quite important. We demonstrated a scenario in which the combined constraints imposed on the Attacker determine the engagement outcome and make capture possible despite the Attacker's kinematic advantage.   

\subsection{Free Attacker Terminal Speed Solutions}
\label{sec:FATSSolutions}

Relaxing the ATS constraint enables the Attacker a richer variety of evasion maneuvers. We discovered four types of such solutions (labeled Types A to D), which we describe here:

\subsubsection{Type A Game Solution}
\label{sec:TypeASolutions}
\begin{figure}[bt!]
\centering    
\subfigure[Players trajectories.]{\label{fig:typeA_trajs}\includegraphics[width=0.45\textwidth]{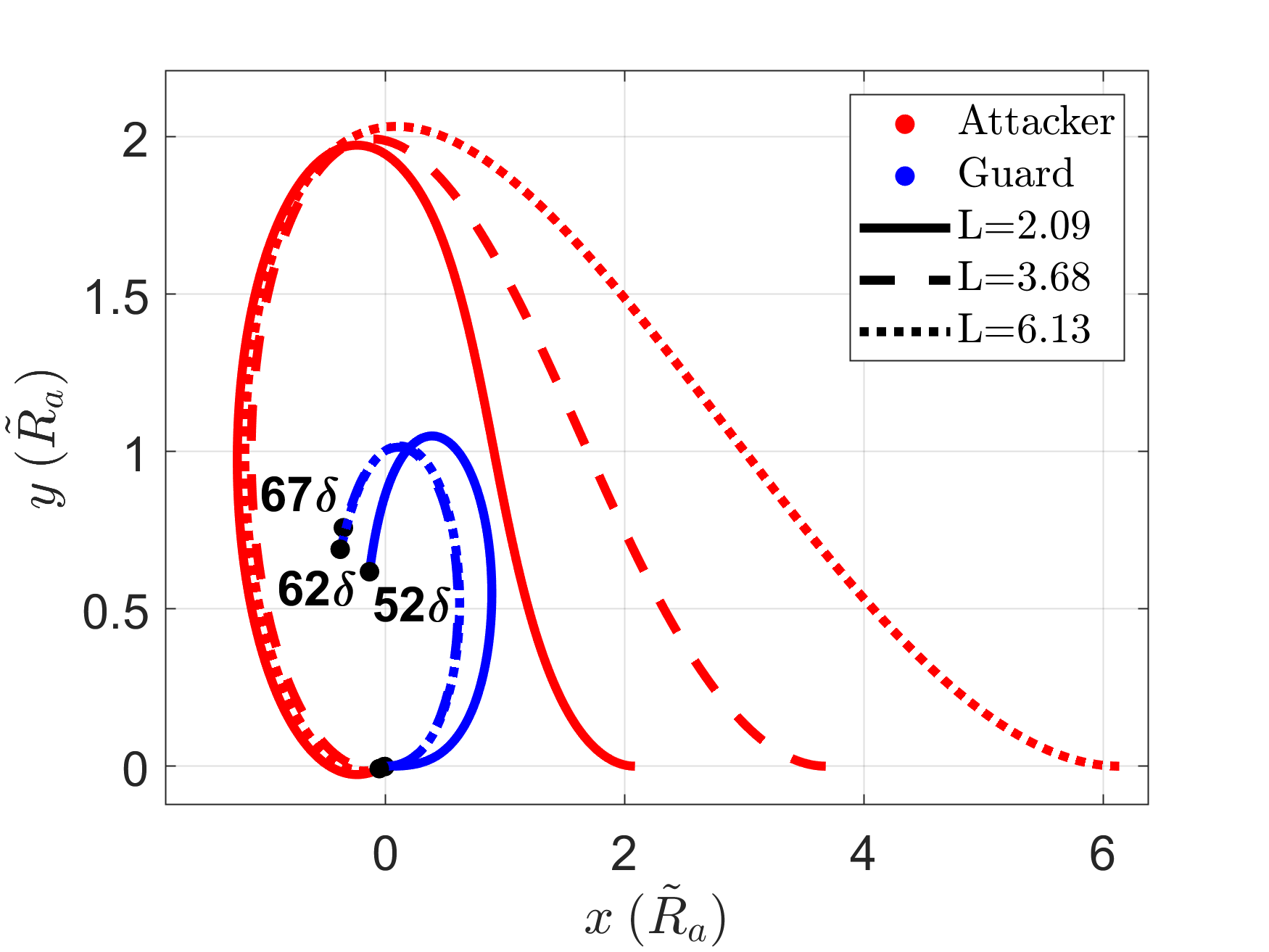}}
\subfigure[Players speed.]{\label{fig:typeA_speed}\includegraphics[width=0.45\textwidth]{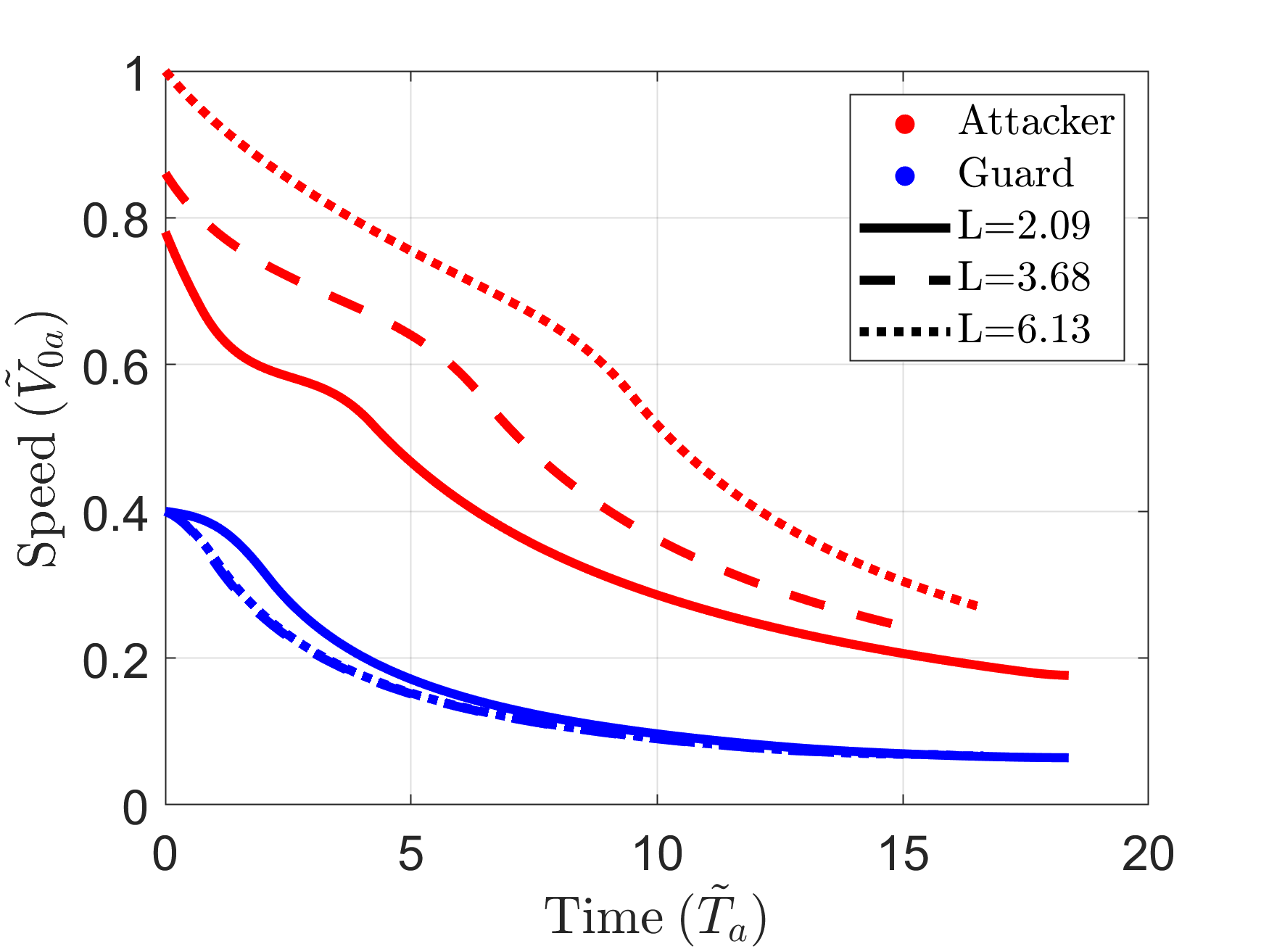}}
\subfigure[Flyby variable.]{\label{fig:typeA_mu}\includegraphics[width=0.45\textwidth]{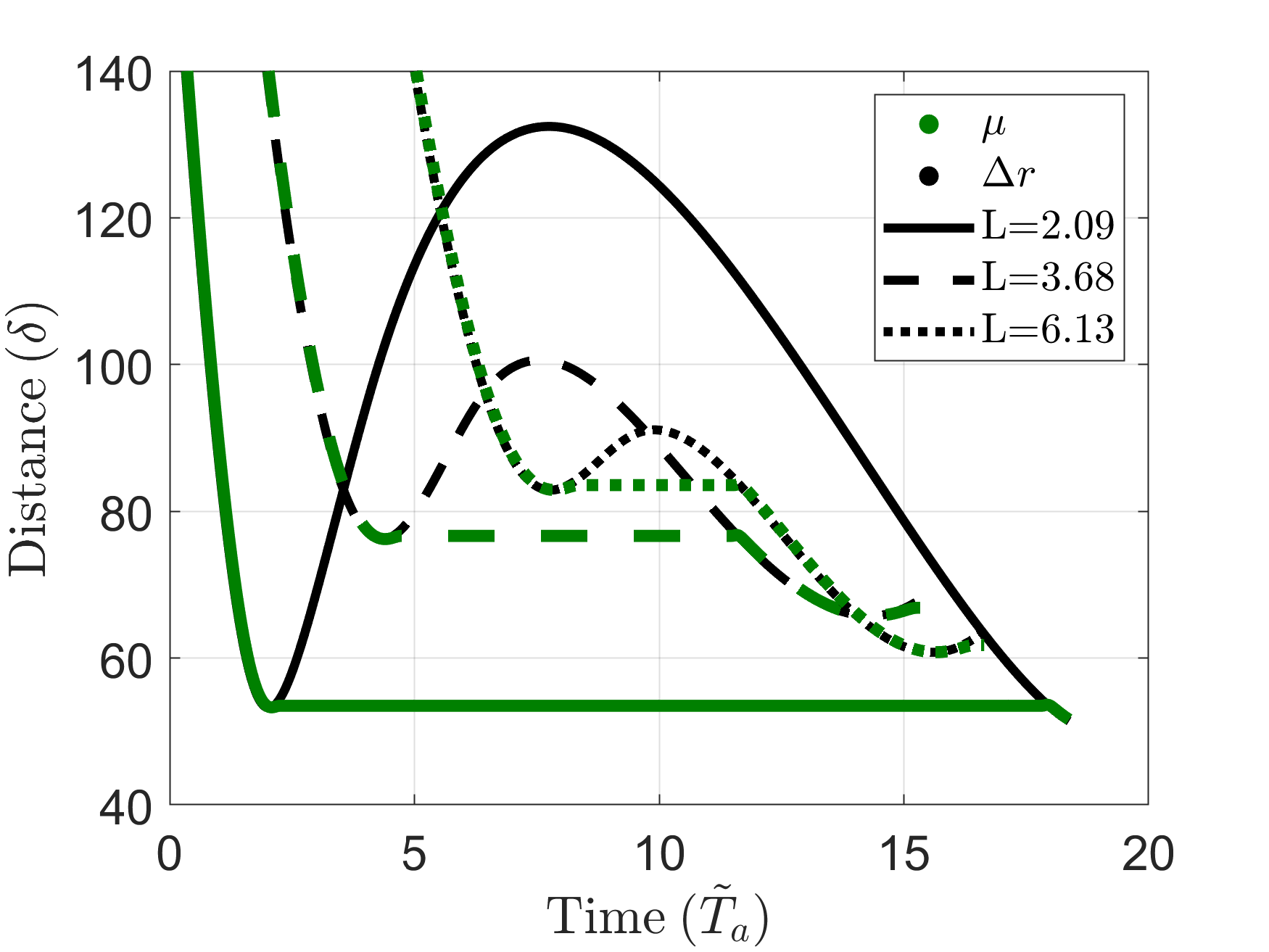}}
\subfigure[Players Controls.]{\label{fig:typeA_control}\includegraphics[width=0.45\textwidth]{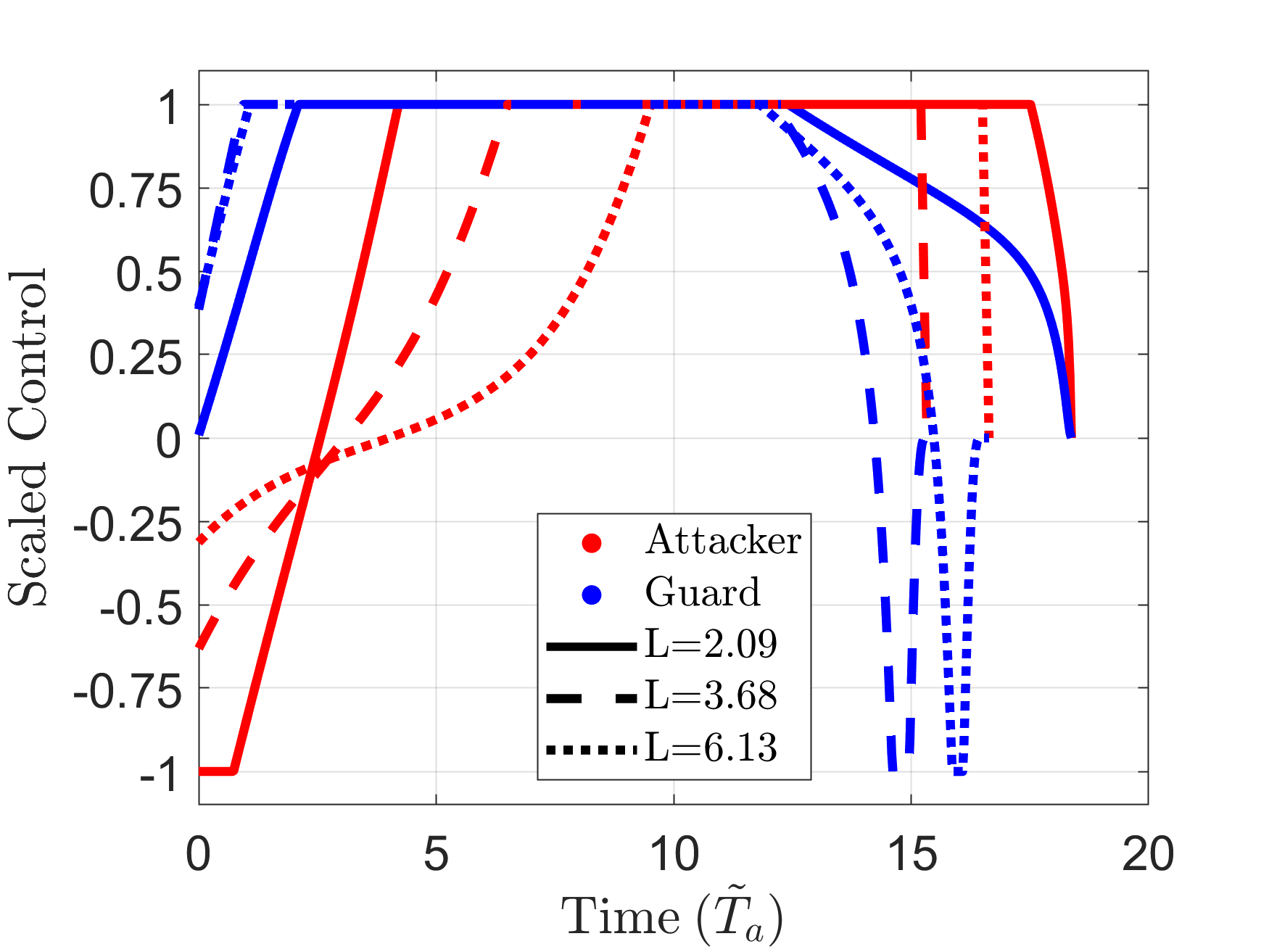}}
\caption{Type A solutions.}
\label{fig:typeA_solutions}
\end{figure}

Type A game solutions are presented in \figref{fig:typeA_solutions}. Figure \ref{fig:typeA_trajs} presents the players' trajectories for game solutions with different launch ranges $L$. The Attacker is again plotted in red, and the Guard is in blue. Black markers on the players' trajectories mark the flyby points, and the terminal flyby variable $\mu(t_f)$ is denoted in $\delta$'s. The Attacker turns initially to the right, flies almost straight for a long distance, and then makes a maximal left turn towards the Target. The Guard tries to follow with a left turn but misses the Attacker by tens of $\delta$'s. The flyby time is just before the Attacker arrives at the Target.  

Type A solutions are similar to the long Bang-Zero-Bang (BZB) solutions discussed in Appendix \ref{sec:AppendixBZB}. BZB solutions are Attacker and Guard trajectories resulting from a degenerate problem in which the Attacker is allowed two maximal-turn flight segments, separated by a straight-flight segment. The Guard is only allowed to execute a constant (possibly maximal) turn. With these constraints, the Attacker's only degrees of freedom are the turn directions and the duration of the first turn. The Guard's only degree of freedom is its constant control. Going back to Eq. \ref{eq:optimal_controls3} we note that the region $\left(-\bar{\varphi}_i,\bar{\varphi}_i\right)$ of moderated control vanishes in the small $\Cd{i}$ limit. In this low induced drag to lift ratio limit, only zero or maximal controls are optimal, and BZB solutions (with maximal Guard turns) thus solve FBDP as well. We used this fact to provide the CNLP solver with effective initial guess solutions in the low induced drag limit. Away from this limit, we used an iterative process that started with a low induced drag solution and gradually adjusted the problem parameters, using the previous iteration solution as a guess for the current one.

Additional information for the Type A solutions can be found in the other three frames of \figref{fig:typeA_solutions}. Figure \ref{fig:typeA_speed} shows the players' speeds, which decrease considerably with flight time. Figure \ref{fig:typeA_mu} shows the flyby distance variable $\mu(t)$, compared with the actual players' separation $\Delta r(t)$. Note that these solutions have two flyby points: the primary flyby point that determines the game value just before termination and an earlier point with a larger flyby distance. In spite of this complexity, we again observe that $\mu$ records the minimal separation as a function of time with good accuracy ($\approx 2\delta$). The difference between the flyby distances at the two points becomes smaller with launch range and equates at $L_A \approx 2$. Type A solutions do not exist below this launch range because the Guard has no reason to continue the chase past the first flyby point once it becomes the primary, and transitions to a speed conservation strategy after the first flyby. 

Figure \ref{fig:typeA_control} shows the players' scaled controls. The Attacker's control is moderated as the launch range increases because it has more time to execute its evasion maneuver. Both players again terminate with zero controls, in agreement with the transversality boundary conditions in Eq. \ref{eq:transversality_bc} and the optimal controls Eq. \ref{eq:optimal_controls3}. 

\subsubsection{Type B Game Solutions}
\label{sec:TypeBSolutions}
\begin{figure}[bt!]
\centering    
\subfigure[Players trajectories.]{\label{fig:typeB_trajs}\includegraphics[width=0.45\textwidth]{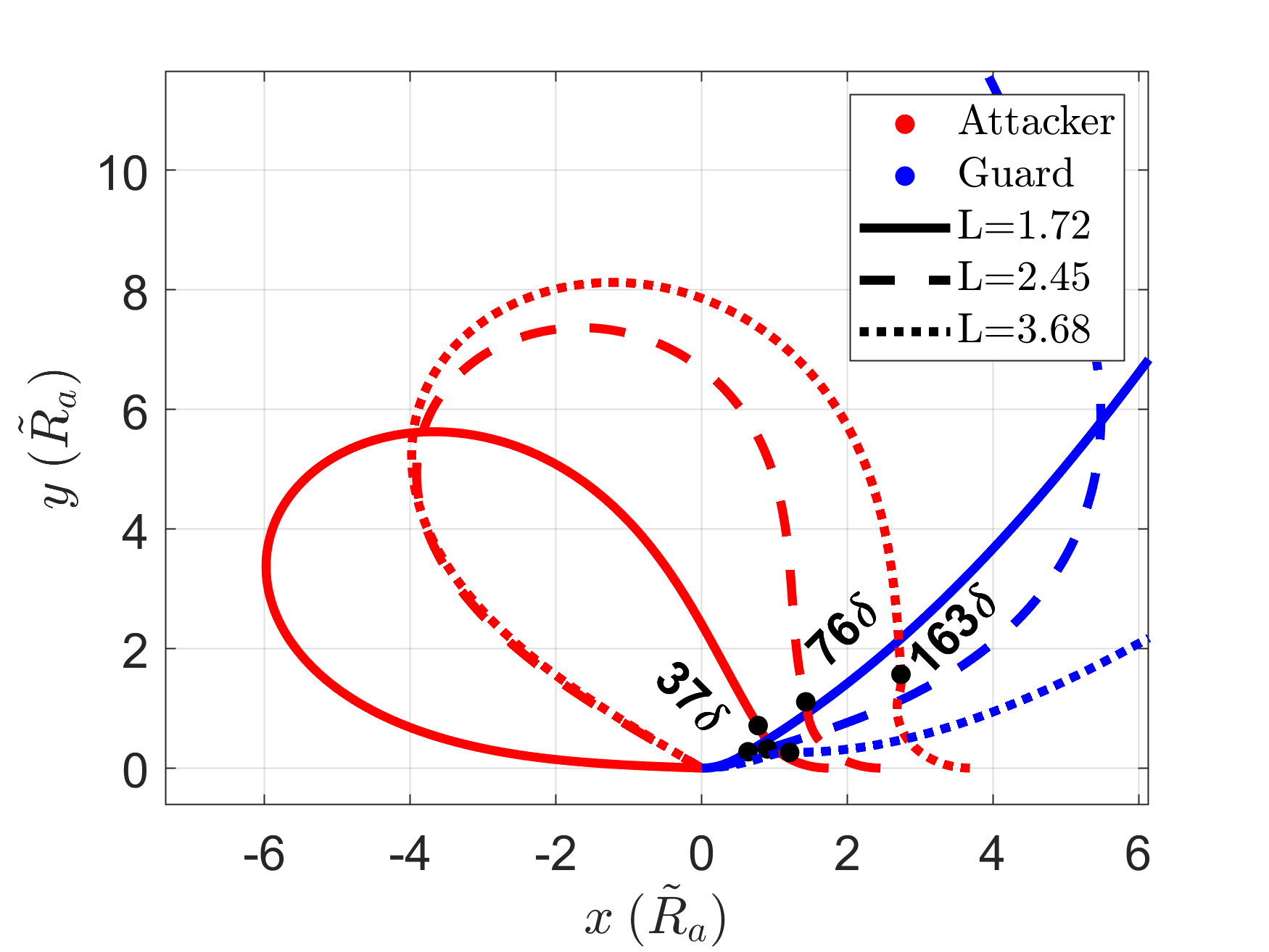}}
\subfigure[Players speeds.]{\label{fig:typeB_speed}\includegraphics[width=0.45\textwidth]{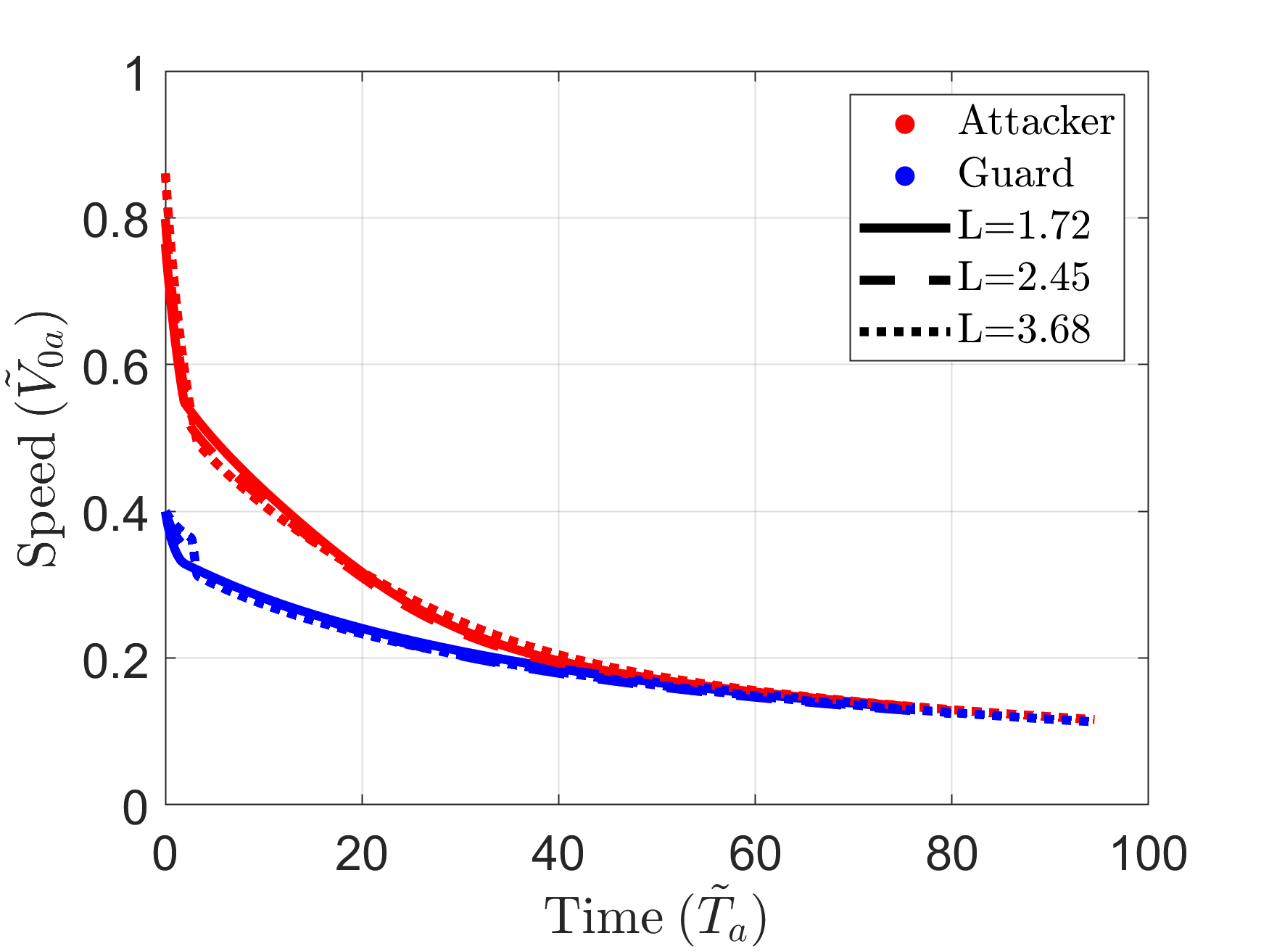}}
\subfigure[Flyby variable.]{\label{fig:typeB_mu}\includegraphics[width=0.45\textwidth]{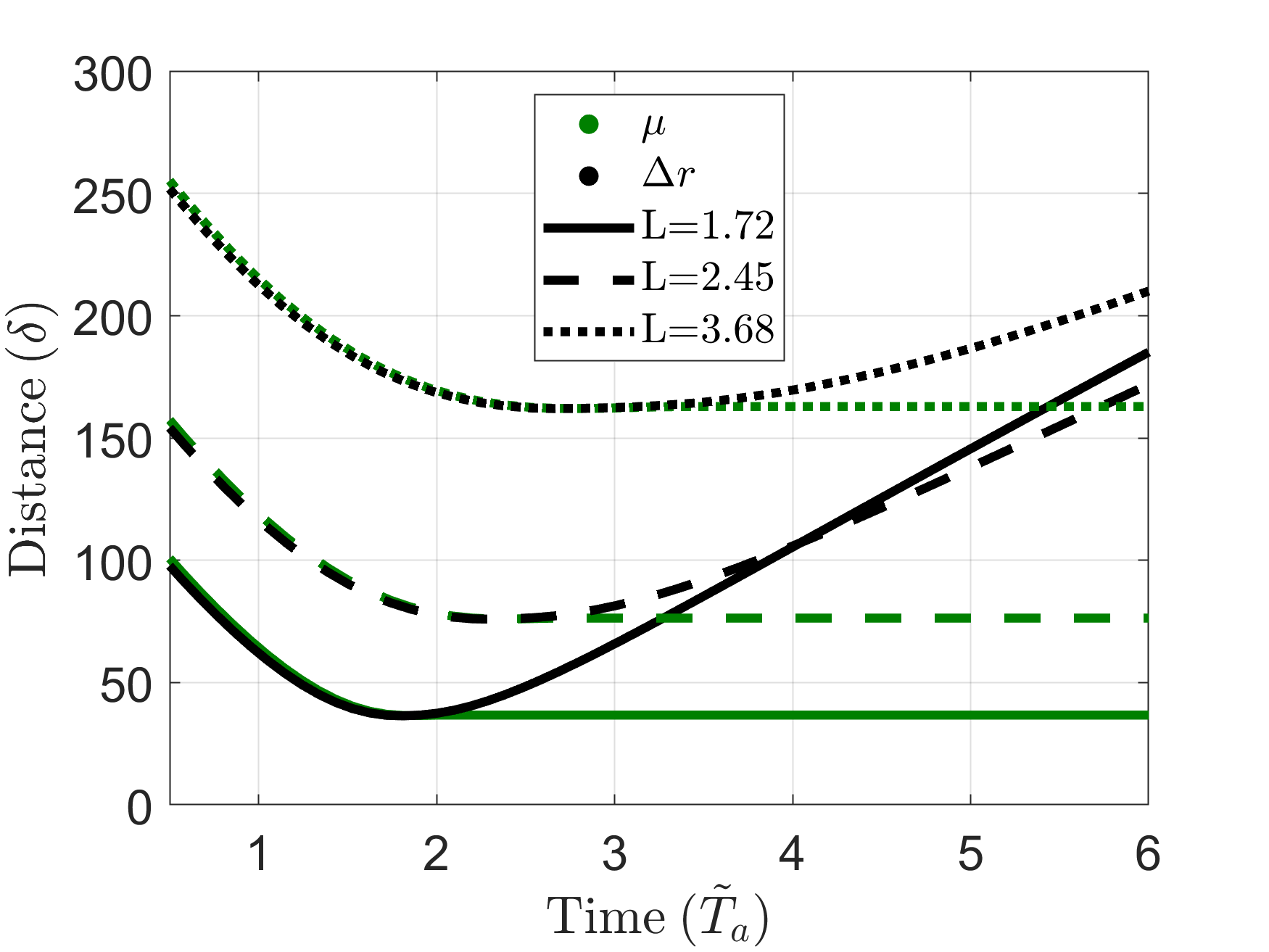}}
\subfigure[Players Controls.]{\label{fig:typeB_control}\includegraphics[width=0.45\textwidth]{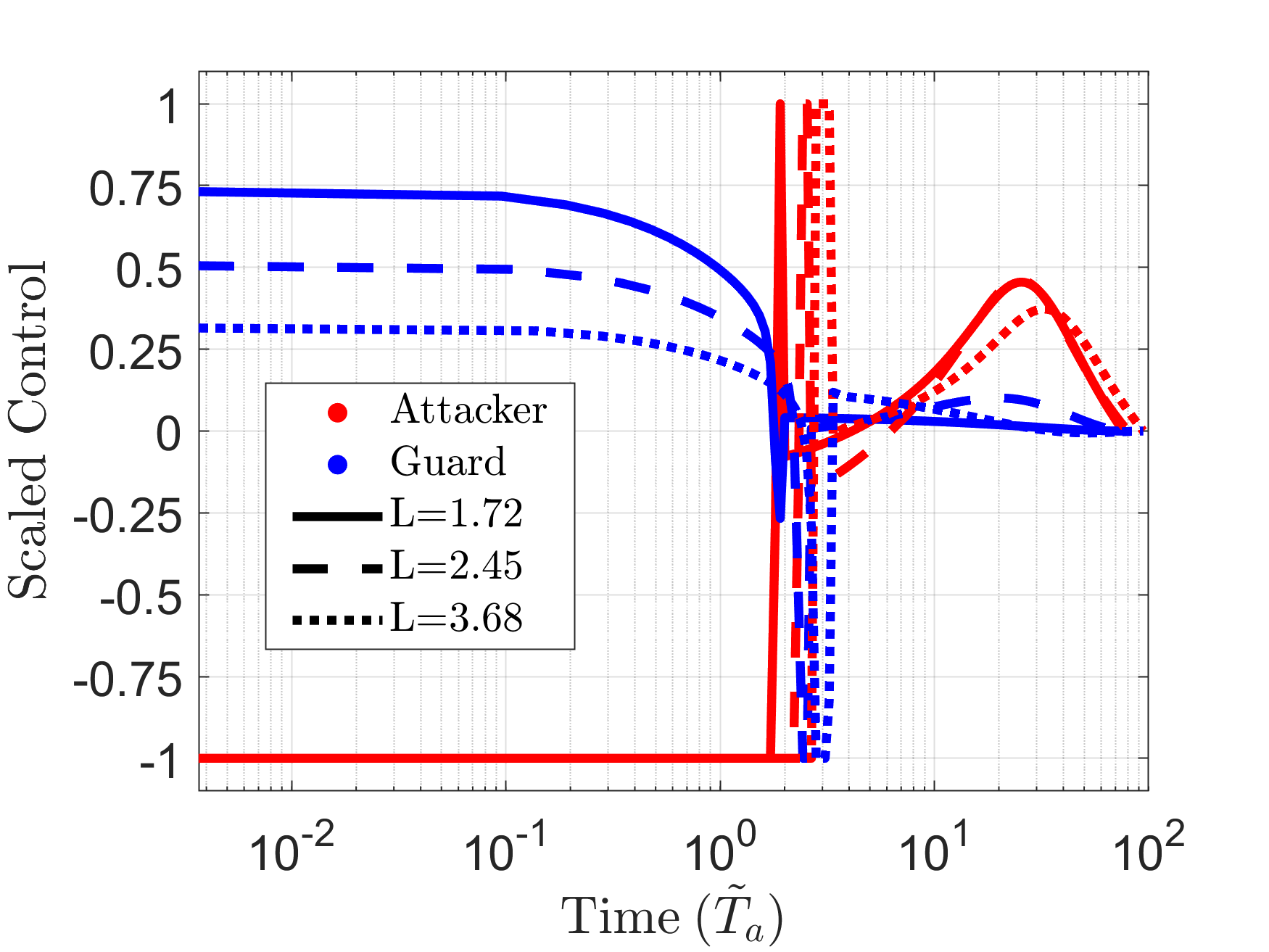}}
\caption{Type B solutions.}
\label{fig:typeB_solutions}
\end{figure}

In type B solutions, the Guard aims for the first flyby point. Examples are presented in \figref{fig:typeB_solutions}. The Attacker executes a maximal right turn up to the flyby point, then makes a long left turn that brings him to the Target (\figref{fig:typeB_trajs}). The Guard makes a moderate left turn towards the flyby point and then continues straight to conserve speed. Flyby distances are many tens of $\delta$'s, flight times are approaching 100, and final speeds of the two players drop below 0.2 (\figref{fig:typeB_speed}). 
The players scaled controls shown in \figref{fig:typeB_control} are plotted against a logarithmic time axis to better capture their behavior. The flyby variable $\mu$ shown in \figref{fig:typeB_mu} is focused on the flyby time. The examples in \figref{fig:typeB_solutions} were selected to match the examples of  \figref{fig:typeA_trajs} and to show an example of a type B solution with $L<L_A$. Guard's post-flyby behavior in these examples exhibits some artifacts due to our limited numerical accuracy at long flight times. These artifacts do not significantly affect the flyby distances.

\begin{figure}[bt!]
\centering    
\subfigure[Flyby distance.]{\label{fig:typesAandBcomp_flybydist}\includegraphics[width=0.45\textwidth]{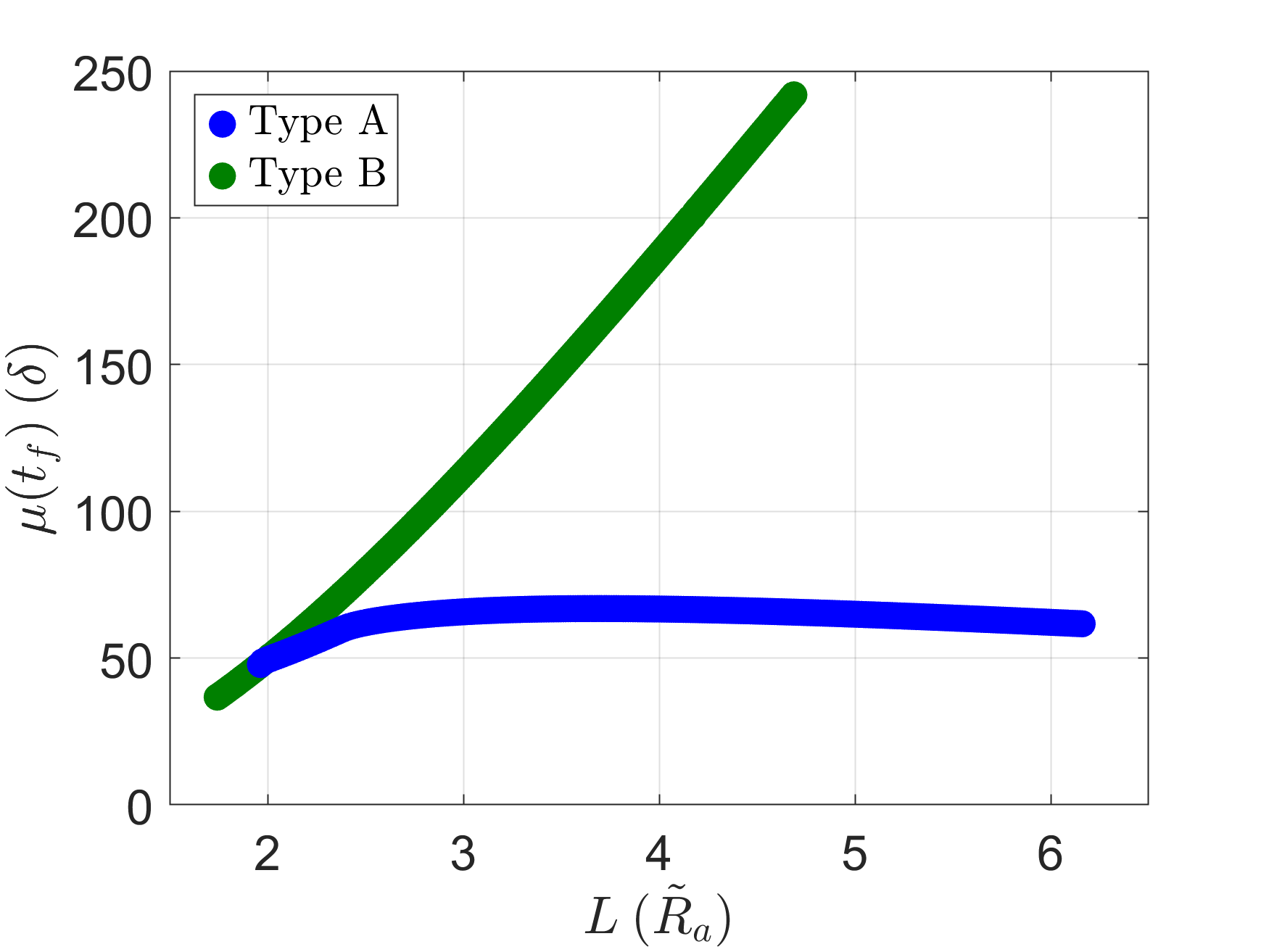}}
\subfigure[Attacker final speed.]{\label{fig:typesAandBcomp_finalvel}\includegraphics[width=0.45\textwidth]{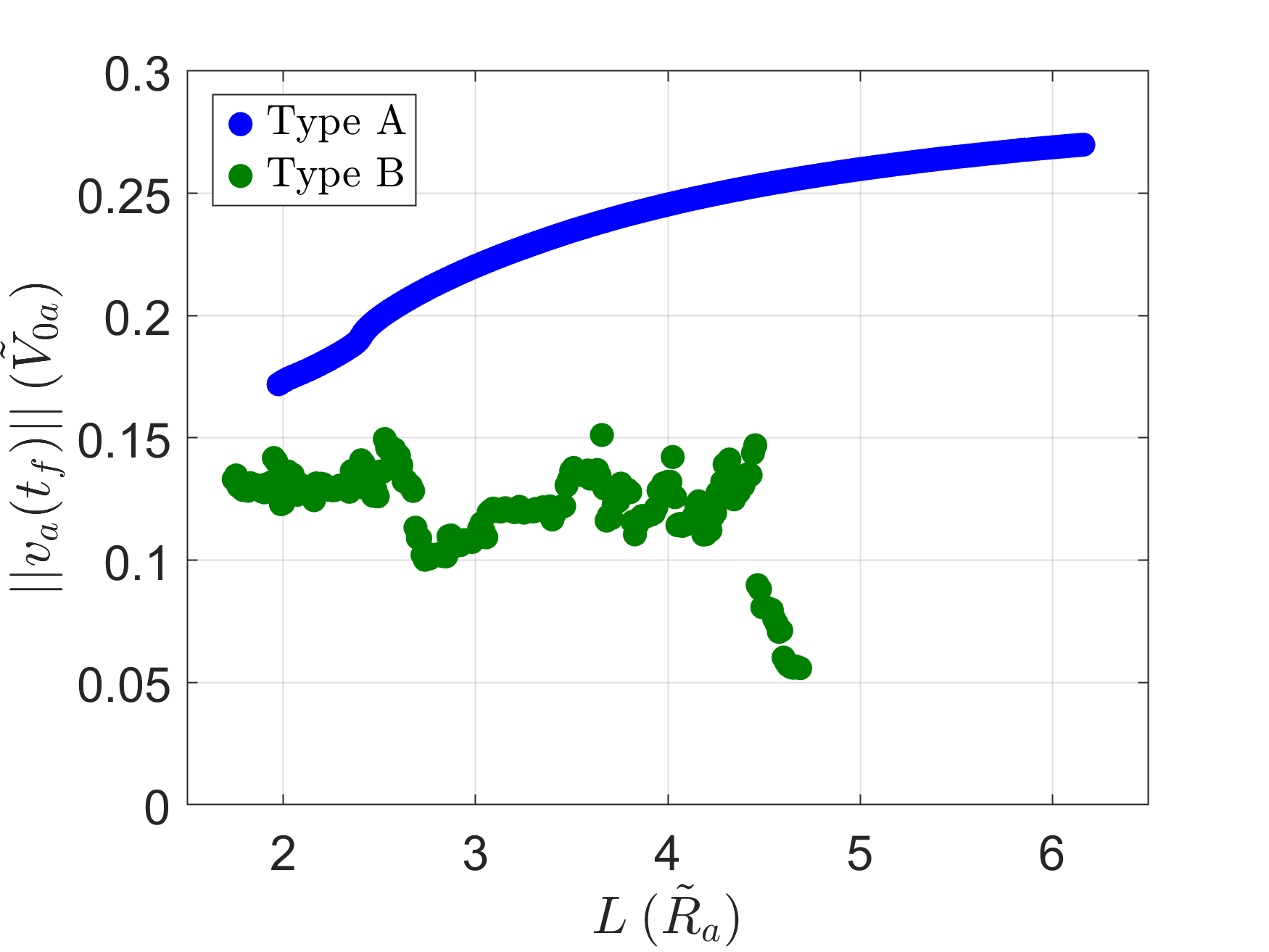}}
\caption{Types A and B engagement performance.}
\label{fig:typeAandBcompare}
\end{figure}

Figure \ref{fig:typeAandBcompare} compares the performance of types A and B in terms of flyby distances and final Attacker speeds. Below $L_A\approx2$, only type B exists. At $L_A$ the two solution types provide the same flyby distance, and above $L_A$, the two types coexist. Type A offers the Attacker flyby distances of 50-70 $\delta$ with terminal speeds of 0.17-0.27. 
Type B offers larger flyby distances and slower terminal speeds. We speculate that Type A can be enforced by the Guard by continuing his turn past the first flyby point. The alternative is that the Attacker can enforce type B by escaping to a great distance and then returning. But this could have also been done in type A, turning type A solutions into something different. The existence of type A therefore suggests that an escape maneuver is not beneficial and that the selection mechanism is the Guard's decision to continue to the second flyby point.    

\subsubsection{Type C Game Solutions}
\label{sec:TypeCSolutions}
\begin{figure}[bt!]
\centering    
\subfigure[Players trajectories.]{\label{fig:typeC_trajs}\includegraphics[width=0.45\textwidth]{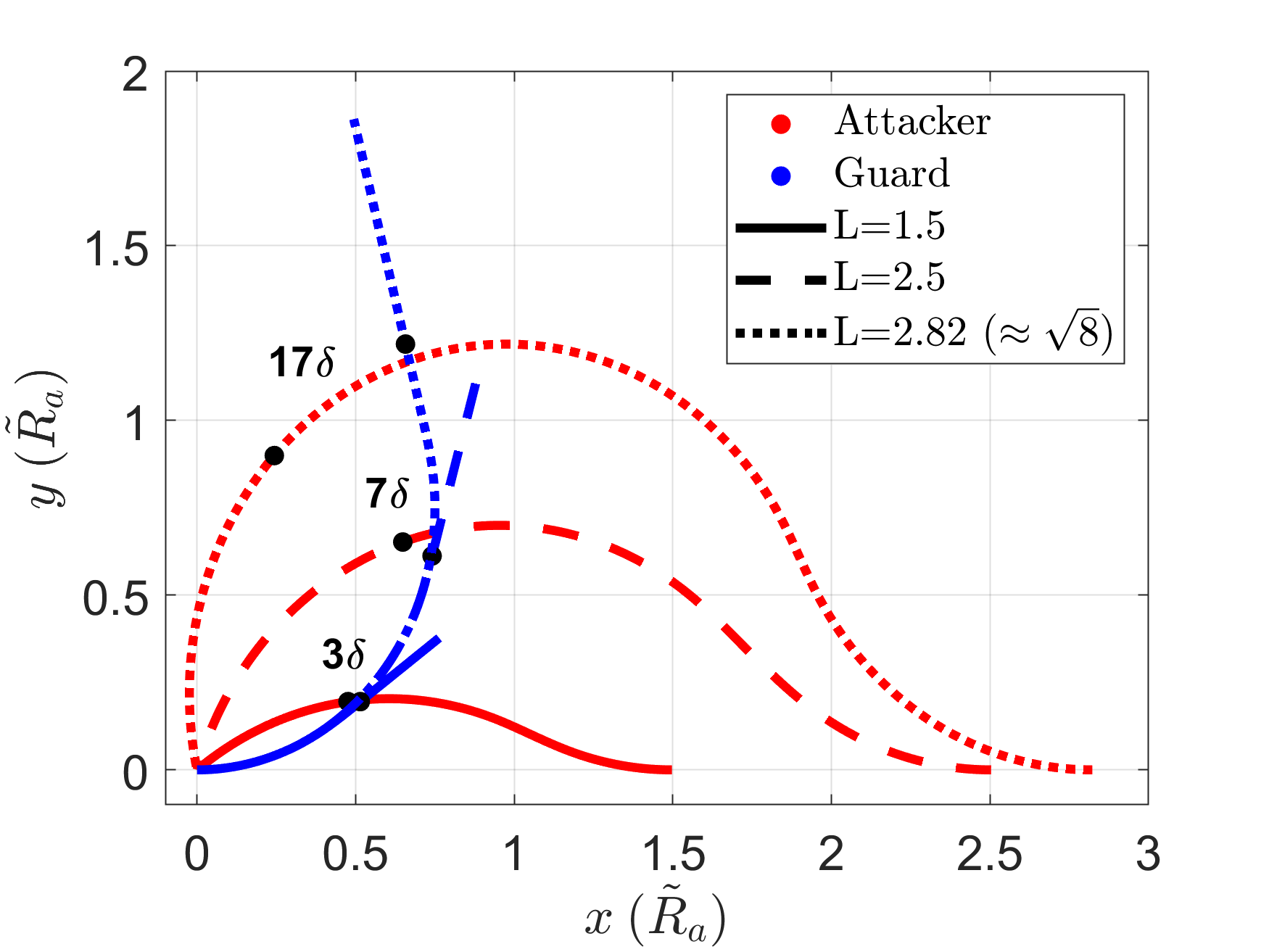}}
\subfigure[Players speed.]{\label{fig:typeC_speed}\includegraphics[width=0.45\textwidth]{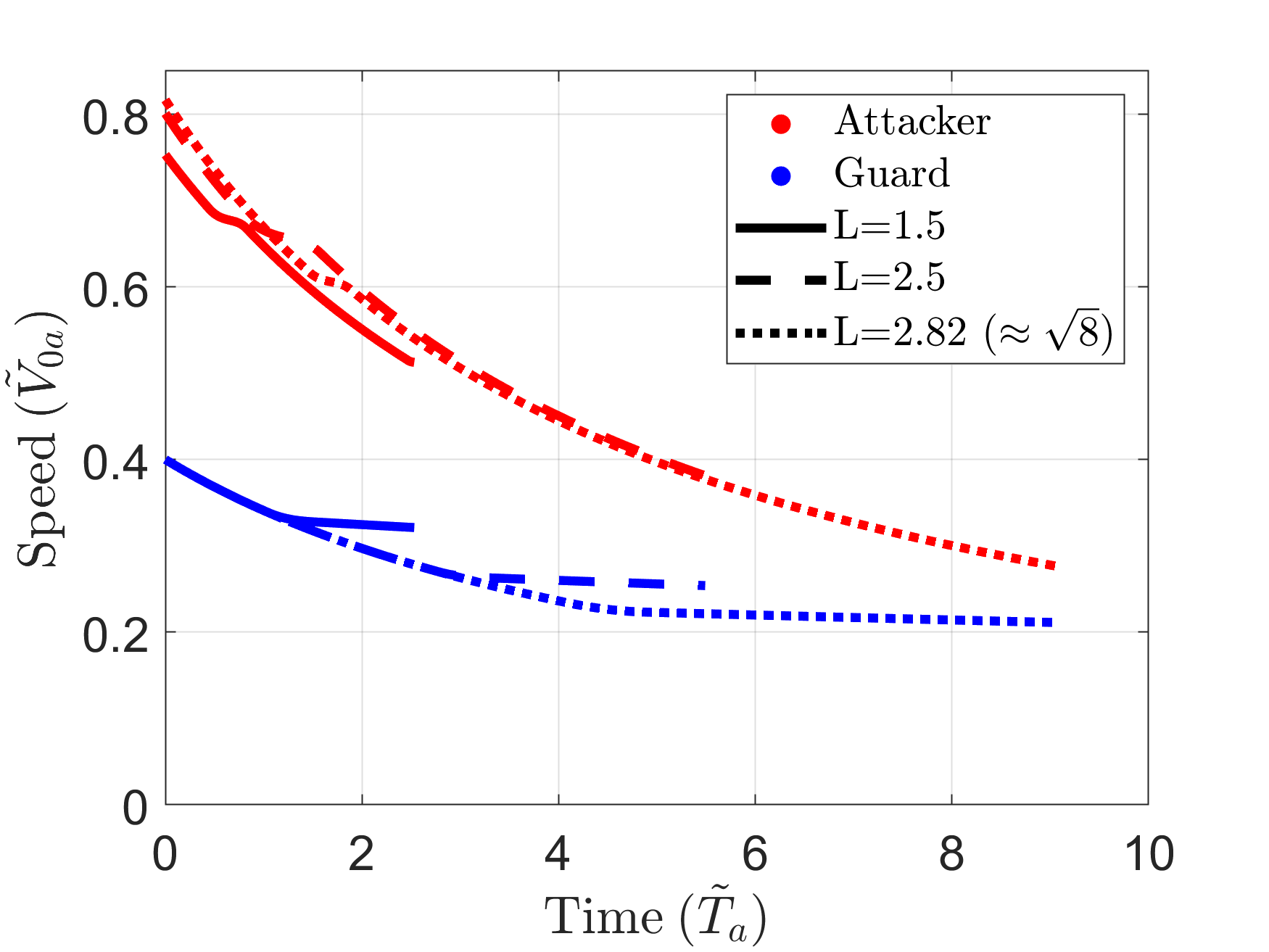}}
\subfigure[Players Controls.]{\label{fig:typeC_control}\includegraphics[width=0.45\textwidth]{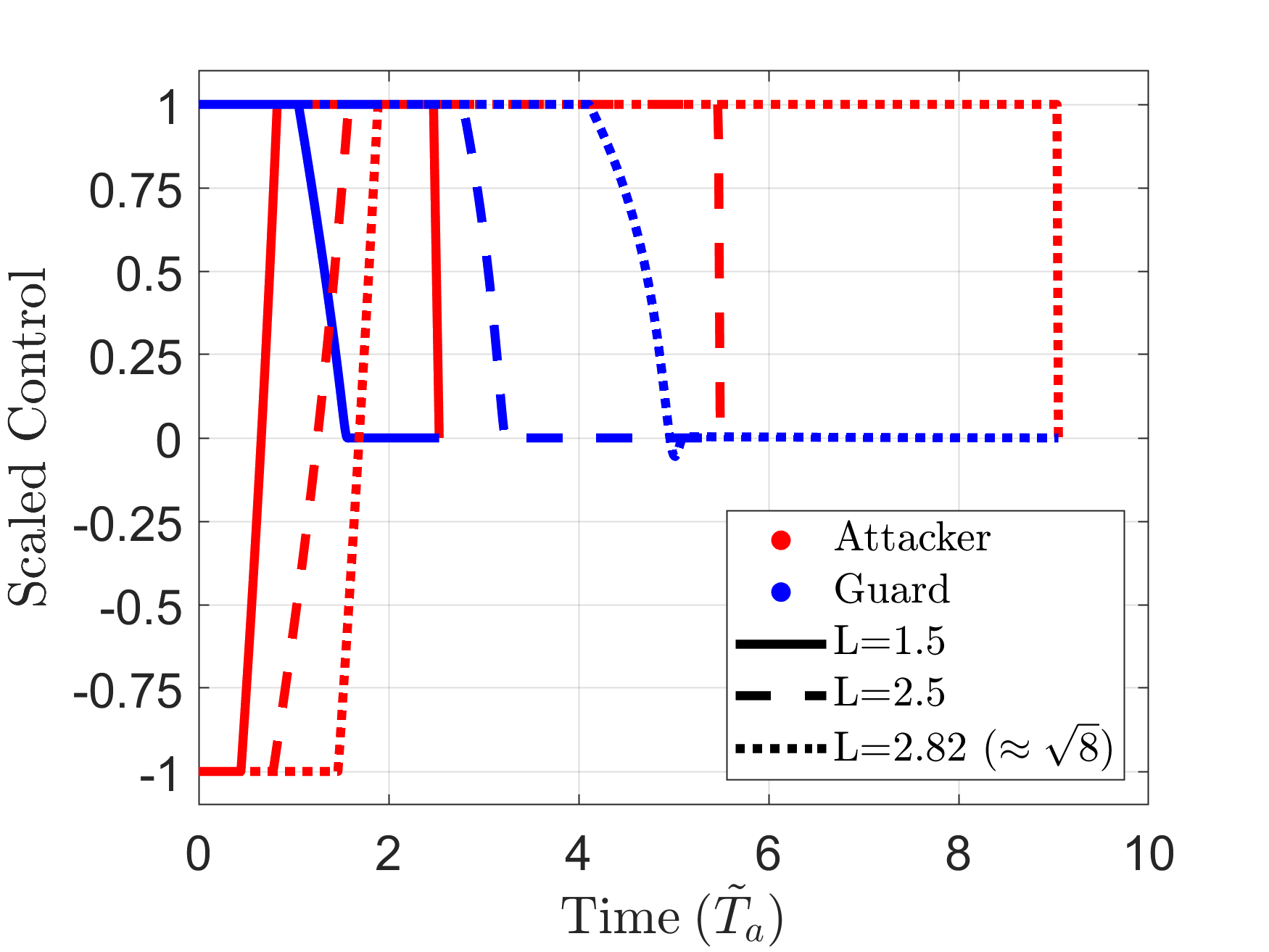}}
\subfigure[Flyby variable.]{\label{fig:typeC_mu}\includegraphics[width=0.45\textwidth]{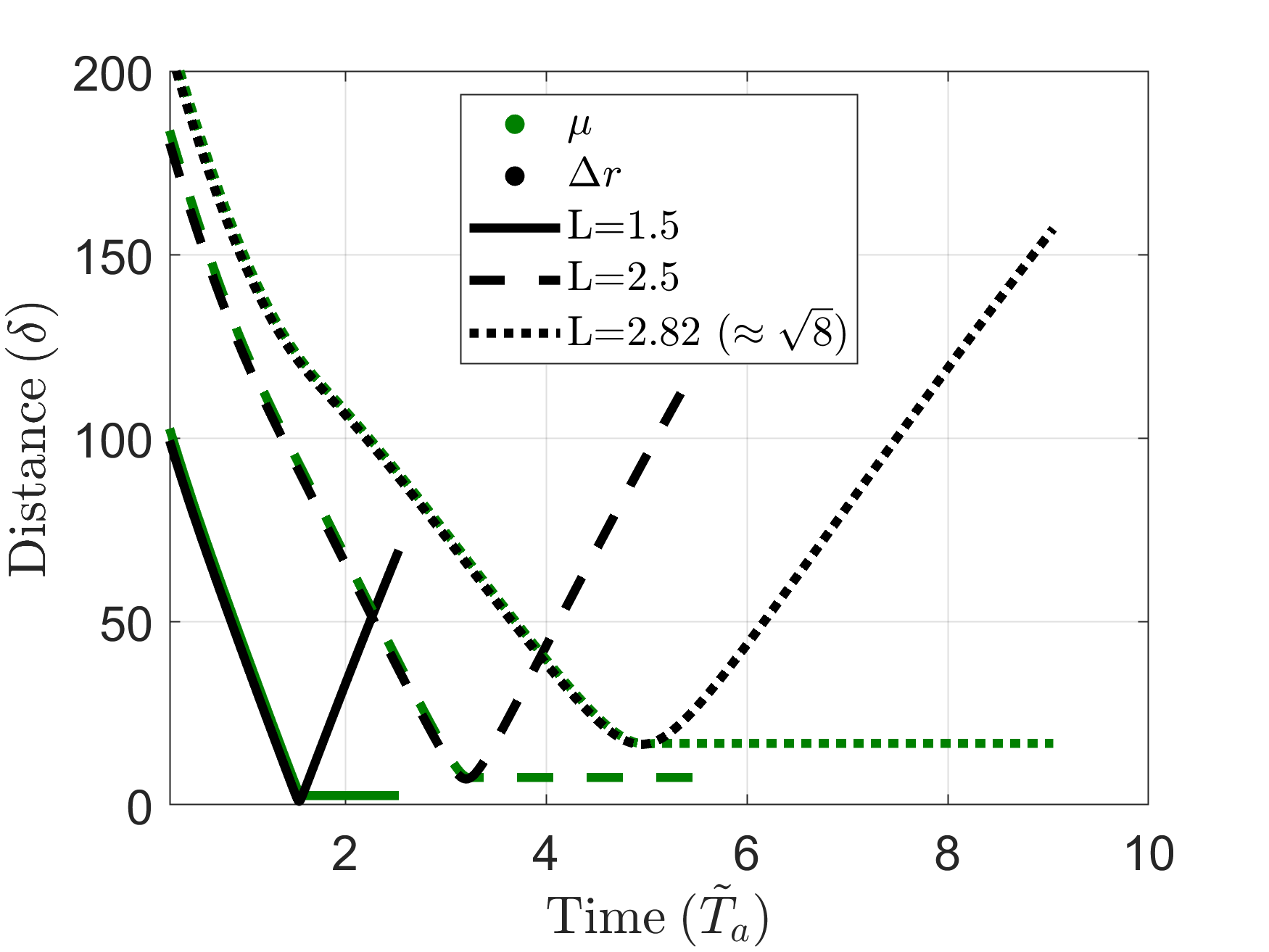}}
\caption{Type C solutions.}
\label{fig:typeCsol}
\end{figure}

Figure \ref{fig:typeCsol} presents the third type of solutions that we found. In order to avoid complete capture, we reduced the Guard's AoA limiter $\umaxi{g}$ from $30\degree$ to $20\degree$. The Guard's scaled parameters are thus modified to $\Cdz{g}=0.046$, $\Cd{g}=0.279$.

Notice that the controls in \figref{fig:typeC_control} are nearly BZB. Type C solutions are indeed almost identical to the short BZB trajectories of Appendix \ref{sec:AppendixBZB}. Moreover, in agreement with the BZB analysis, type C solutions do not exist for launch ranges $L\geq\sqrt{8}$. As explained in Appendix \ref{sec:AppendixBZB}, at this launch range, the short BZB branch merges with the long BZB and is no longer a saddle point.  

\begin{figure}[bt!]
\centering    
\subfigure[Flyby distance.]{\label{fig:typeC_performance_mu}\includegraphics[width=0.45\textwidth]{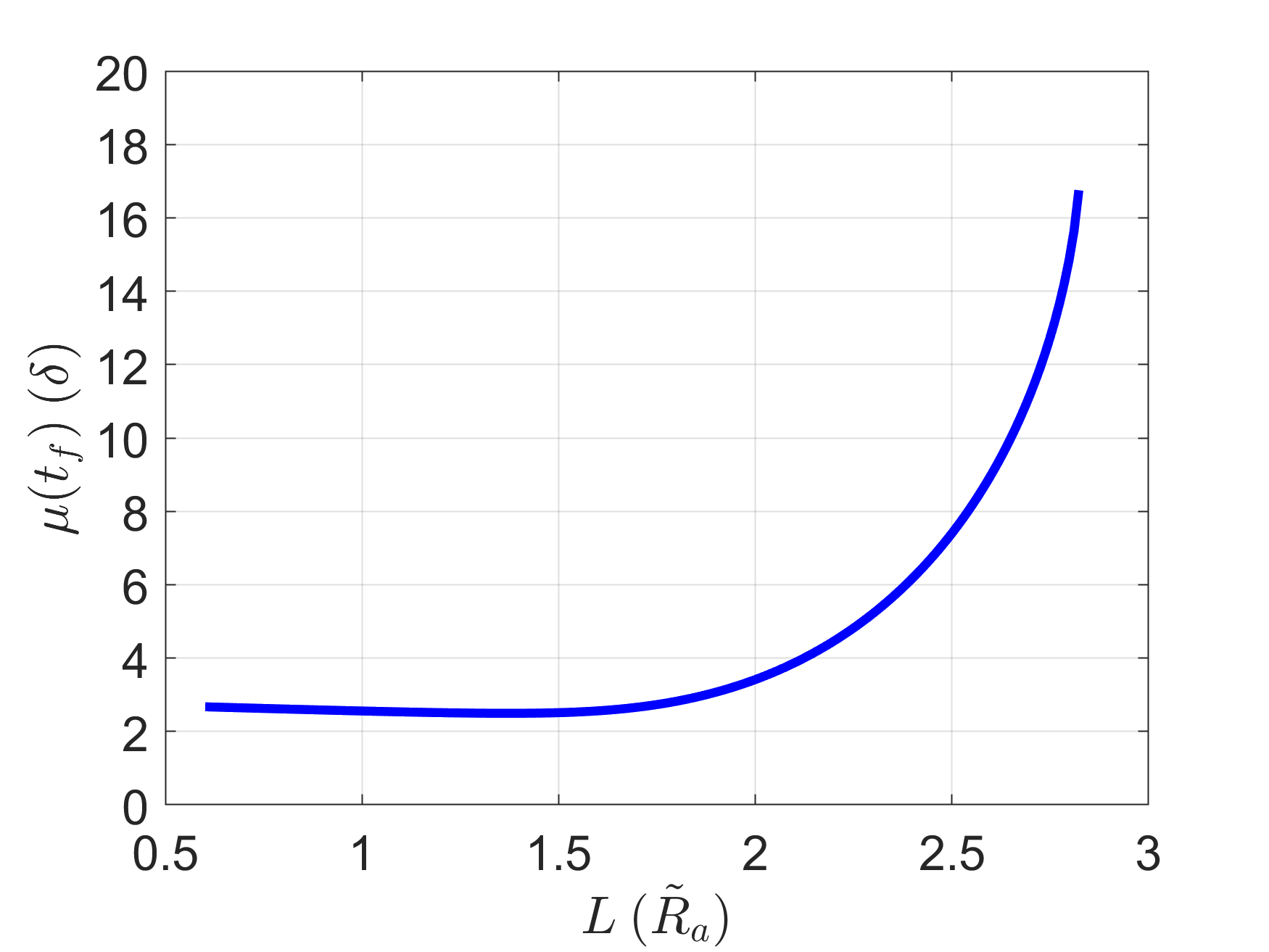}}
\subfigure[Final Attacker speed.]{\label{fig:typeC_performance_speed}\includegraphics[width=0.45\textwidth]{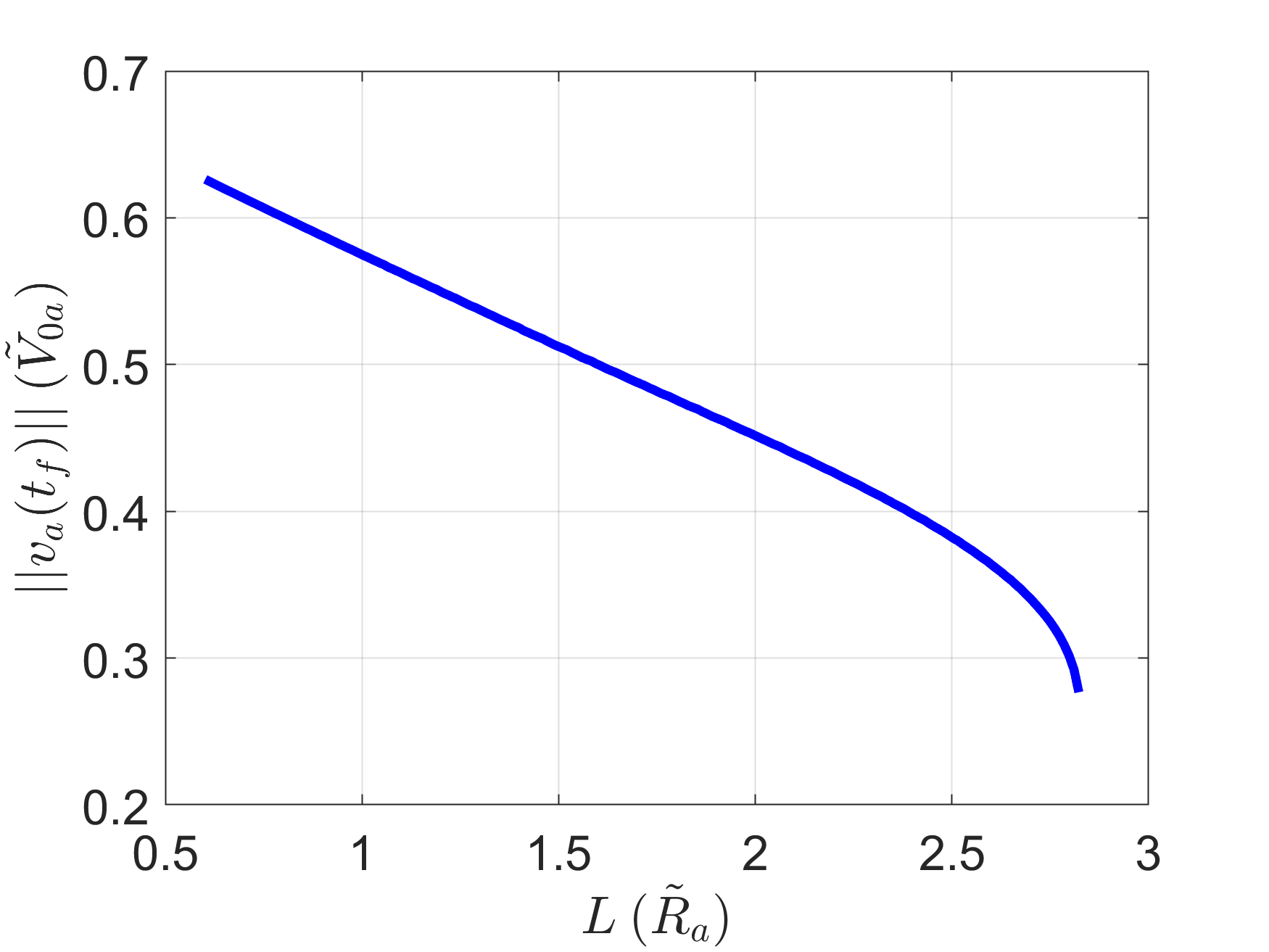}}
\caption{Type C solutions performance.}
\label{fig:typeC_performance}
\end{figure}

Figure \ref{fig:typeC_performance} presents the performance of type C solutions in terms of flyby distances and final Attacker speeds. Even with the reduced Guard AoA limiter, type C offers the Attacker flyby distances shorter than solution types A or B. For launch ranges below 1.5, he is actually captured. The selection between type C and A/B solutions is done by shortening or prolonging the Attacker's initial right turn. Type C solutions are therefore local saddle points because the Attacker can enforce types A/B.

\subsubsection{Type D Game Solutions}
\label{sec:TypeDSolutions}
\begin{figure}[bt!]
\centering    
\subfigure[Players trajectories.]{\label{fig:typeD_trajs}\includegraphics[width=0.45\textwidth]{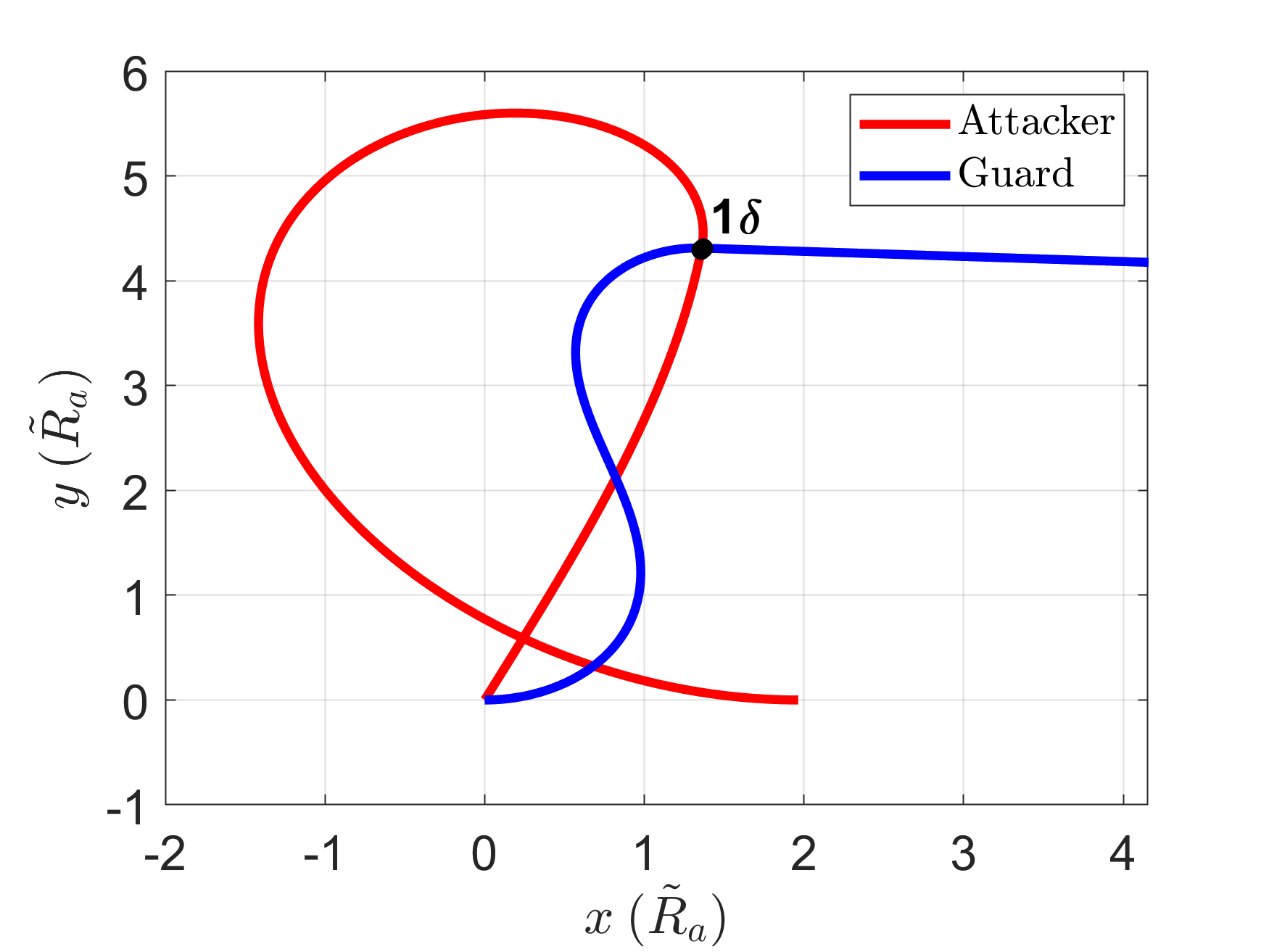}}
\subfigure[Players speed.]{\label{fig:typeD_speed}\includegraphics[width=0.45\textwidth]{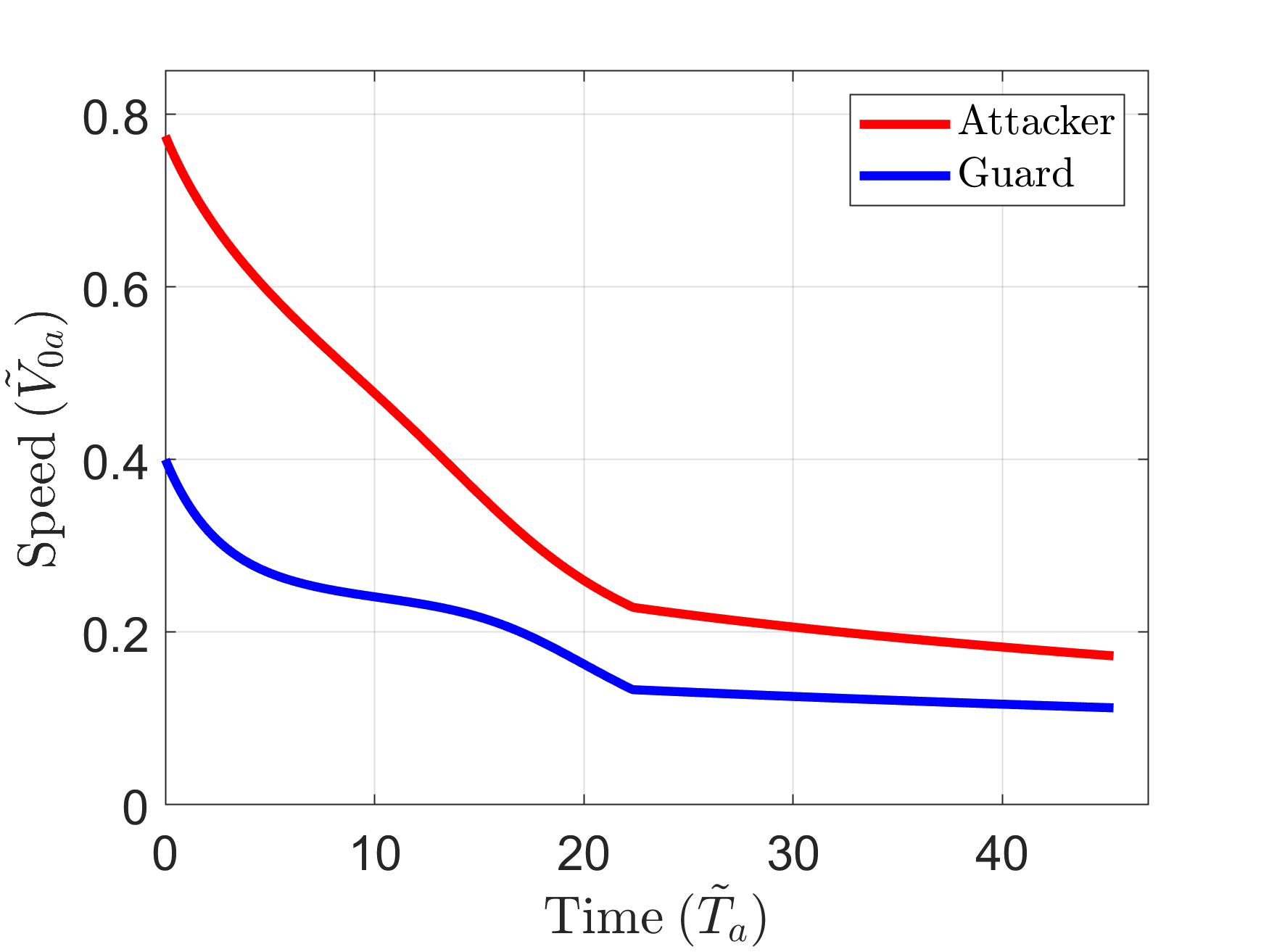}}
\subfigure[Players controls.]{\label{fig:typeD_control}\includegraphics[width=0.45\textwidth]{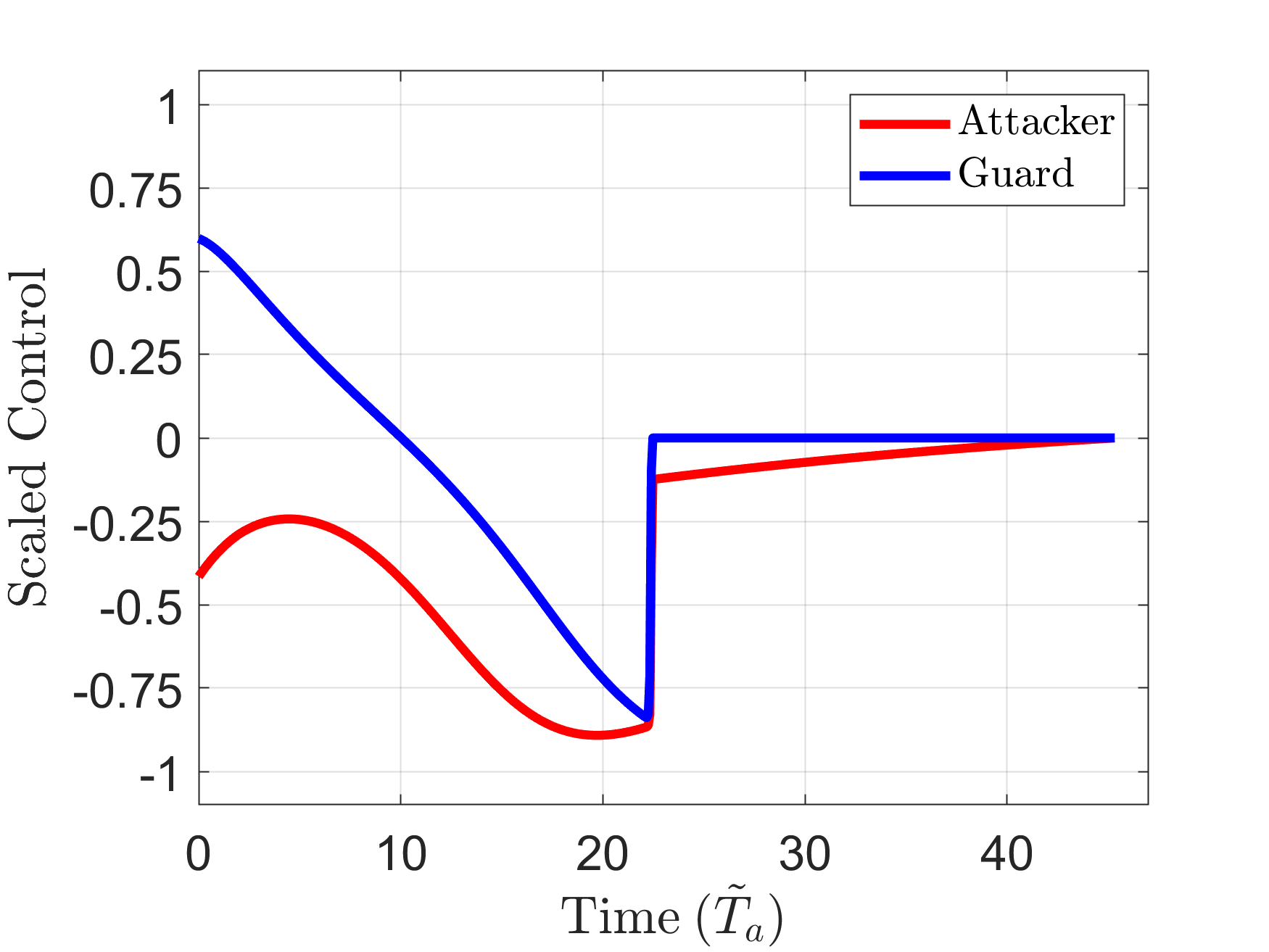}}
\subfigure[$\mu$]{\label{fig:typeD_mu}\includegraphics[width=0.45\textwidth]{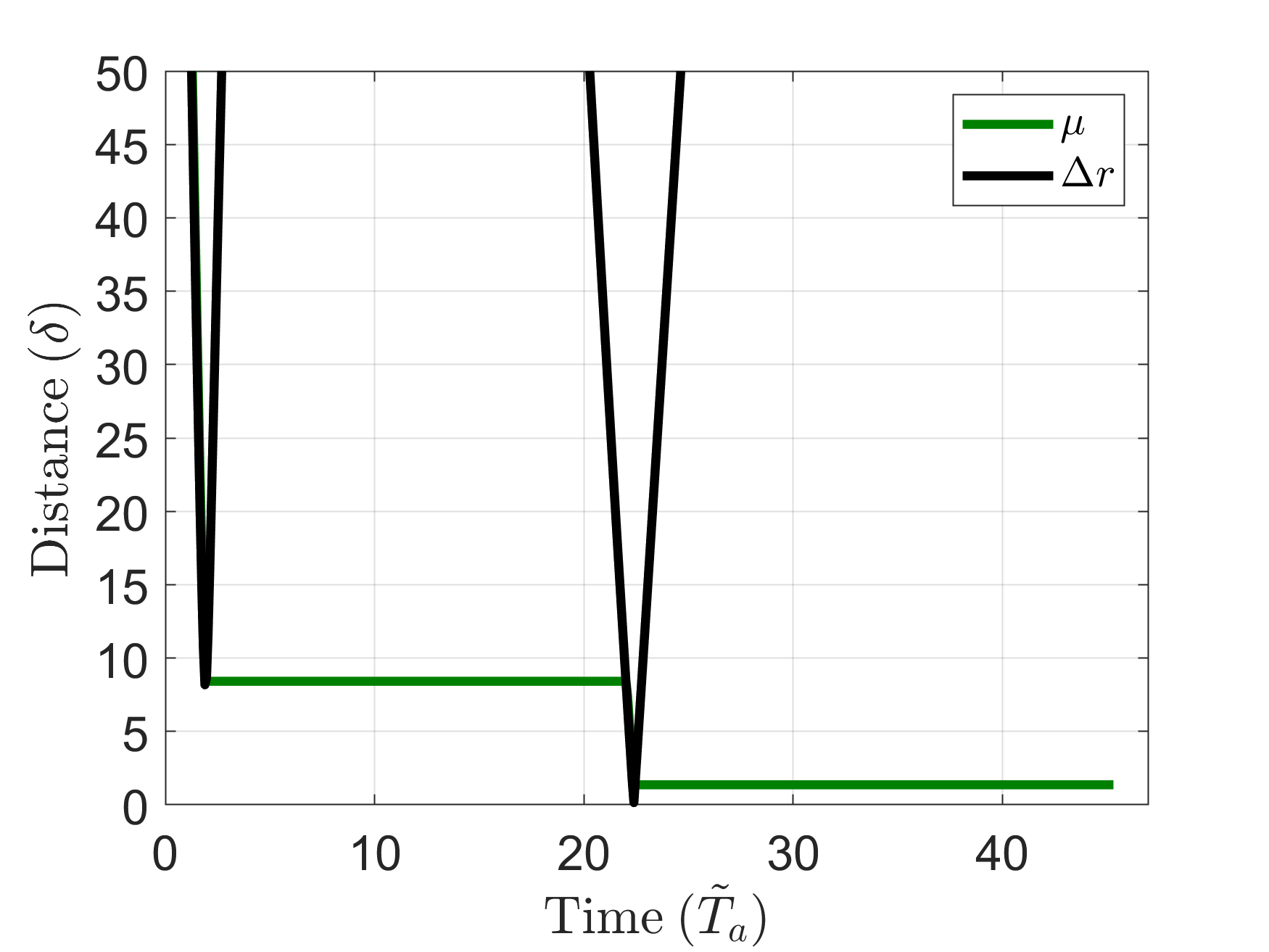}}
\caption{Type D solution.}
\label{fig:typeD_solution}
\end{figure}

The solutions we presented above consist of two basic Attacker maneuvers: evasion followed by Target homing. In contrast, the solution shown in \figref{fig:typeD_solution} consists of a single one. The parameters are identical to those used in \secref{sec:TypeASolutions} except for the Attacker's speed cost weight, which is reduced to $\phi_{v_a}=0.014$.
The Guard launches at the range $L\approx2$. 
It is clearly a local saddle point, since the Attacker is caught and, in the same situation with a larger speed cost ($\phi_{v_a}=0.37$), he could have escaped with a flyby distance of more than 40$\delta$ (see \figref{fig:typeAandBcompare}). 
The Guard here performs a rather elaborate trajectory, which may be interpreted as first blocking a direct approach and then flying to intercept the loop. Since the Attacker here is caught, his strategy is to exhaust the Guard (whose speed cost weight is $\phi_{v_g}=0.37$). The existence of this local saddle point suggests that there may be others. Given the complex behavior exhibited here, it may be difficult to identify and map all such solutions.  

To summarize this section, relaxing the ATS constraint leads to a rich variety of solutions. We identified four distinct types and suspect that there may be additional types. Types A and B offer the Attacker large flyby distances at the cost of slow terminal speeds. Flyby distances in B are larger, and we currently don't know how the selection mechanism between the two types works. Namely, it is not clear if A can be enforced by the Guard (making B a local saddle point) or B can be enforced by the Attacker (making A a local saddle point). We speculate that the former is true, but cannot prove so at this time. Type C is clearly a local saddle point that the Attacker may select but yields suboptimal performance. Type D is another local saddle point that exhibits an elaborate chase that terminates with a capture. These local saddle points attest to the problem complexity, but are impractical as guidance laws. The obtained solutions should therefore be interpreted as local solutions of the nonlinear saddle-point game.

\section{Conclusions}
\label{sec:Conclusions}
Guarding a target against an Attacker was posed as a differential game. The game is played in a scenario in which the Attacker has a kinematic advantage over the Guard. If the Guard's mission was simply defined to catch the Attacker, he would fail. To mount a practical defense, we defined a viable defense mission and tasked the Guard with preventing the Attacker from hitting a specific Target. Namely, the Guard should catch the Attacker only on Target ingress trajectories. The Attacker's mission is to hit the Target without being intercepted. Attacker's trajectories that terminate elsewhere are considered a defense success.   

The above missions were translated into a novel differential game that we termed Flyby Distance Pursuit (FBDP). In contrast to classical Guarding a Target (GaT) problems \cite{isaacs_differential_1999,pachter_differential_2017} that terminate upon capture, FBDP is terminated when the Attacker reaches the target. The game value is determined by the flyby distance, i.e., the minimum separation between players throughout their flight. This value is made available upon game termination by introducing an auxiliary state variable that records the minimum separation. 

FBDP was formulated for two players moving in a constant-altitude plane.
Playing on a plane is a main simplification in our work, which should obviously be extended in the future to three-dimensional space. The game equations were numerically solved for a reference scenario using the PSOPT optimal control solver. We applied two types of boundary conditions to the Attacker terminal speed: Free and Constrained Attacker Terminal Speeds (FATS and CATS). CATS solutions are practical in the sense that the Attacker is conscious of his terminal speed and maintains his kinematic advantage. Varying the terminal speed constraint, we were able to draw the Attacker's tradeoff between flyby distance and terminal speed. Our results show that reasonable terminal speed constraints lead to capture. Namely, defense is possible despite the Attacker's kinematic advantage. 

FATS solutions are quite impractical due to long flight times and slow terminal speeds. However, they allowed us to identify different possible pursuit-evasion strategies and to demonstrate the problem's complexity. We identified four distinct types of FATS solutions. In contrast to CATS solutions, the Attacker here can evade the Guard and reach the Target with large flyby distances (at the cost of terminal speed).    

\section*{Acknowledgements}
This work was supported by an internal Technion - Israel Institute of Technology fund.

\appendix  %This command ends the counting of sections.
\section{Appendix: Bang-Zero-Bang Flyby Distance Pursuit}
\label{sec:AppendixBZB}
\begin{figure}[bt!]
\centering
\includegraphics[width=0.7\textwidth]{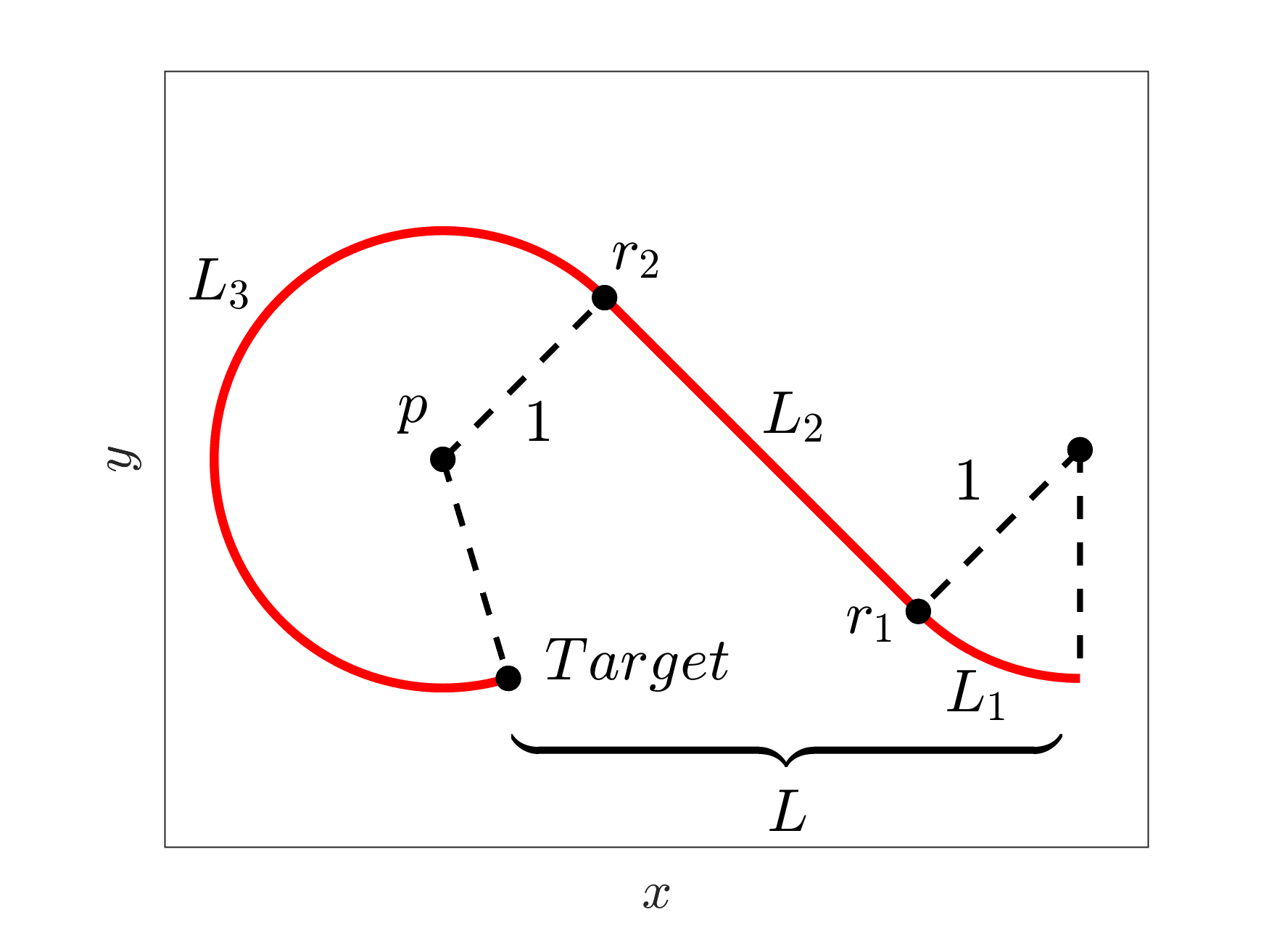}
\caption{BZB Attacker trajectory.}
\label{fig:BZB_illustration}
\end{figure}

It is instructive to study a degenerated version of the pursuit game where the Attacker is allowed a Bang-Zero-Bang (BZB) maneuver consisting of the three flight segments illustrated in \figref{fig:BZB_illustration}: maximal right turn with an arc-length $L_1$, straight flight for a distance $L_2$ and maximal left turn with an arc-length $L_3$. We work in the scaled units of Eq. \ref{eq:scaled_eom} (Attacker turn radius is 1) and consider the same head-on scenario of \secref{sec:FBDPSolutions} where the Guard is co-located with the Target and launches at a range $L$.  

At the end of the first turn, the Attacker is at 
$\ri{1}=\left[L-\sin{L_1},\ \ 1-\cos{L_1}\right]^T$. 
At the beginning of the second turn, he is at
$\ri{2}=\ri{1}-L_2\cdot \left[\cos{L_1},\ \ -\sin{L_1} \right]^T$
and is turning left around the pivot point 
\begin{equation}
\pmb{p}=\ri{2}-\left[ \sin{L_1}, \cos{L_1}\right]^T=\left[L-2 \cdot \sin{L_1}-L_2 \cdot \cos{L_1}, 1-2\cdot\cos{L_1}+L_2\cdot\sin{L_1} \right]^T 
\end{equation}

For the Attacker to reach the Target at the origin, we must have 
$||\pmb{p}||^2=1$.
This is a quadratic equation for $L_2$ with the solutions 
\begin{equation}
\label{eq:L2qe}
L_2=L\cdot\cos{L_1}-\sin{L_1}\pm\sqrt{1-\left(L\cdot\sin{L_1}+\cos{L_1} -2\right)^2}
\end{equation}

\begin{figure}[bt!]
\centering
\subfigure[Attacker BZB straight flight segment $L_2$ as a function of the initial turn length $L_1$.]{\label{fig:L2_vs_L1}\includegraphics[width=0.45\textwidth]{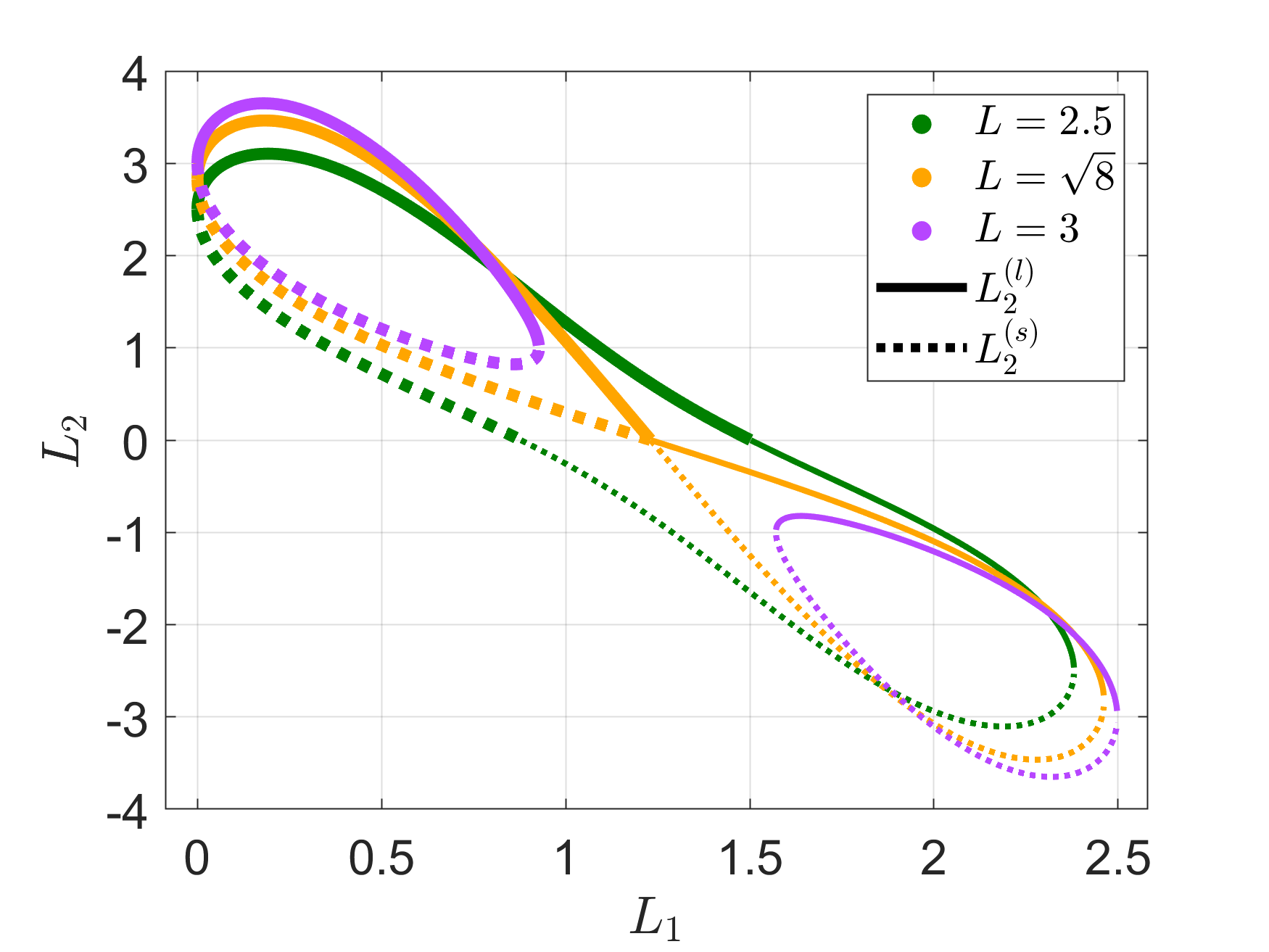}}
\subfigure[Flyby distance as a function of $L_1$.]{\label{fig:miss_vs_L1}\includegraphics[width=0.45\textwidth]{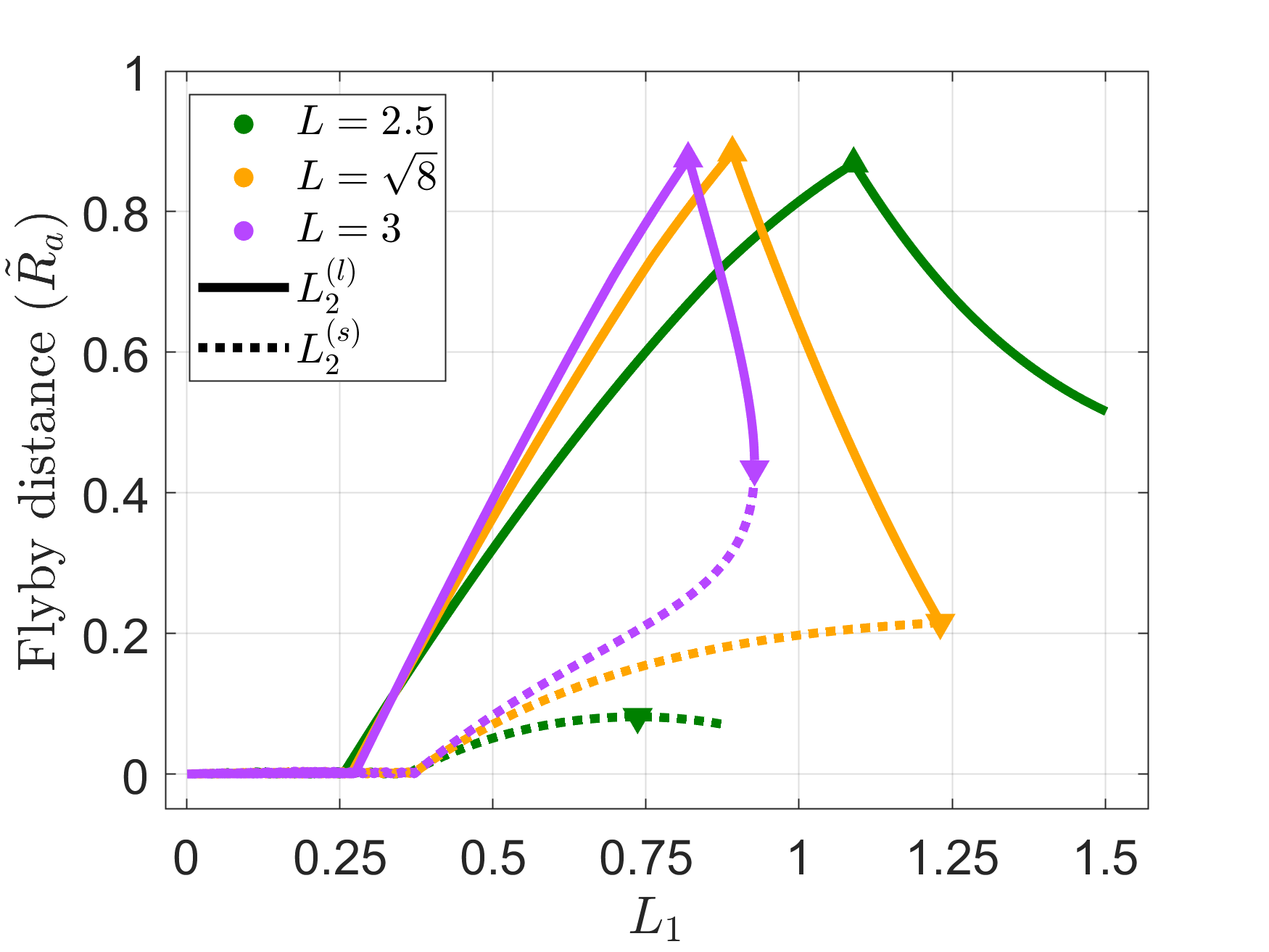}}
\caption{BZB branches.}
\end{figure}

We therefore have two families of Attacker BZB solutions, short ($L_2^{(s)}$) and long ($L_2^{(l)}$), depending on the sign of the discriminant above. Figure \ref{fig:L2_vs_L1} shows the two branches of $L_2$ as a function of $L_1$, for several (scaled) launch ranges $L$. This figure was generated with the same parameters as in \secref{sec:TypeCSolutions}, and we indicate the actual launch ranges in the legend for convenient reference to \figref{fig:typeC_trajs}. Solid/dashed curves show the $L_2^{(l)}$/$L_2^{(s)}$ solutions respectively. Note that we are only interested in solutions with $L_2\geq0$. Negative values imply that the Attacker's velocity is reversed at $\ri{1}$ and are non-physical. We emphasized those positive solutions in \figref{fig:L2_vs_L1} using thicker lines. 

The discriminant in Eq. \ref{eq:L2qe} becomes zero when   
\begin{equation}
\label{eq:zerodet}
L\cdot\sin{{L_1}_{\pm}}+\cos{{L_1}_{\pm}} -2=\pm 1
\end{equation}
The second (minus) option is solved by ${L_1}_-=0$ and ${L_1}_-={L_1}_{max}=atan2(2L,1-L^2)$, the two values that bound the $L_1$ interval with real $L_2$ solutions. Note that $L_2(L_1=0)=L$, which implies straight flight. $L_2({L_1}_{max})=-L$ which is non-physical. 

When $L<\sqrt{8}$ the discriminant is positive for $0<L_1<{L_1}_{max}$ and the two branches  $L_2^{(l)}, \ L_2^{(s)}$ do not meet. At $L=\sqrt{8}$ the two branches coincide with $L_2^{(l)}=L_2^{(s)}=0$ and for $L>\sqrt{8}$ they coincide with  $L_2\neq0$ values. The critical value $L=\sqrt{8}$ may be derived from the ${L_1}_+$ solution of Eq. \ref{eq:zerodet} which we rewrite as
\begin{equation} 
\label{eq:criticalL_1}
L=\frac{3-\cos{{L_1}_+}}{\sin{{L_1}_+}}
\end{equation}
This function has a single minimum in ${L_1}_+\in(0,{L_1}_{max})$, satisfying $\cos{{L_1}_+}=\frac{1}{3}$, which leads to $L=\sqrt{8}$. 

The Guard is allowed a turn with a constant control $\ui{g}\leq1$, which he optimizes for each selection of an Attacker BZB trajectory. Since both players fly with piecewise constant turns, their equations of motion can be integrated analytically. Using Eq. \ref{eq:scaled_eom}, the speed $v_i(t) \equiv ||\vi{i}(t)||$ during a turn with constant $\ui{i}$ obeys 
\begin{equation}
\frac{dv_i(t) }{dt}=-\zeta_i\cdot v_i^2(t) \cdot C_{d_{total}}, \quad C_{d_{total}}=\Cdz{i}+\Cd{i}\cdot \ui{i}^2
\end{equation}
and is solved by
\begin{equation}
v_i(t) =\left[\zeta_i \cdot C_{d_{total}} \cdot (t-t_0)+v_i^{-1}(t_0)\right]^{-1}    
\end{equation}
The traveled distances grow as
\begin{equation}
\label{eq:bzb_traveled_distance}    
l_i(t)=l_i(t_0)+\frac{1}{\zeta_i\cdot C_{d_{total}}}\cdot \left[\log\left(\zeta_i \cdot C_{d_{total}}\cdot(t-t_0)+v_i^{-1}(t_0)\right)+\log v_i(t_0)\right]
\end{equation}

Inverting Eq. \ref{eq:bzb_traveled_distance}, we get the flight time as a function of the traveled distance

\begin{equation}
\label{eq:bzb_flight_time}    
t(l_i)=t_0+\frac{1}{\zeta_i\cdot C_{d_{total}}\cdot v_i(t_0)}\cdot \left\{\exp{\left[\zeta_i \cdot C_{d_{total}}\cdot \left(l_i-l_i(t_0)\right)\right]}-1\right\}
\end{equation}
which allows us to draw the players' trajectories as a function of time and compute flyby distances.  

Under the above restrictions, the Guard's optimal trajectory is found by adjusting $\ui{g}$ to minimize the flyby distance for every choice of an Attacker trajectory. The results are presented in \figref{fig:miss_vs_L1}. The different curves show the Guard-minimized flyby distances as a function of the initial arc-length $L_1$. Solid/dashed curves again denote  ${L}_2^{(l)}$/${L}_2^{(s)}$ solutions. Curves maxima are marked by triangles where $\Delta  / \nabla$ are the maxima of ${L}_2^{(l)}$/${L}_2^{(s)}$ curves, respectively. The maxima of ${L}_2^{(l)}$ solutions are the best result the Attacker can achieve from an initial distance ${L}$. The $\Delta$'s are thus global Guard-Attacker (min-max) saddle points. The green $\nabla$ on the ${L}_2^{(s)}$ curve is a local saddle point. It provides the Attacker with shorter flyby distances than the green $\Delta$, but any variation of $L_1$ at this point will degrade his result.

In contrast, the orange $\nabla$ is a maximum of the orange ${L}_2^{(s)}$ curve but is not a saddle point. It is located on the long ${L}_2^{(l)}$ curve as well, and the Attacker may move from it to the orange $\Delta$ while continually increasing his game value. The same argument holds for the purple curve as well, and for every launch range above $L=\sqrt{8}$. The short BZB saddle point therefore disappears at $L\geq\sqrt{8}$.
    
\begin{figure}[bt!]
\centering    
\subfigure[$L=2.5$.]{\label{fig:bzb_traj1}\includegraphics[width=0.3\textwidth]{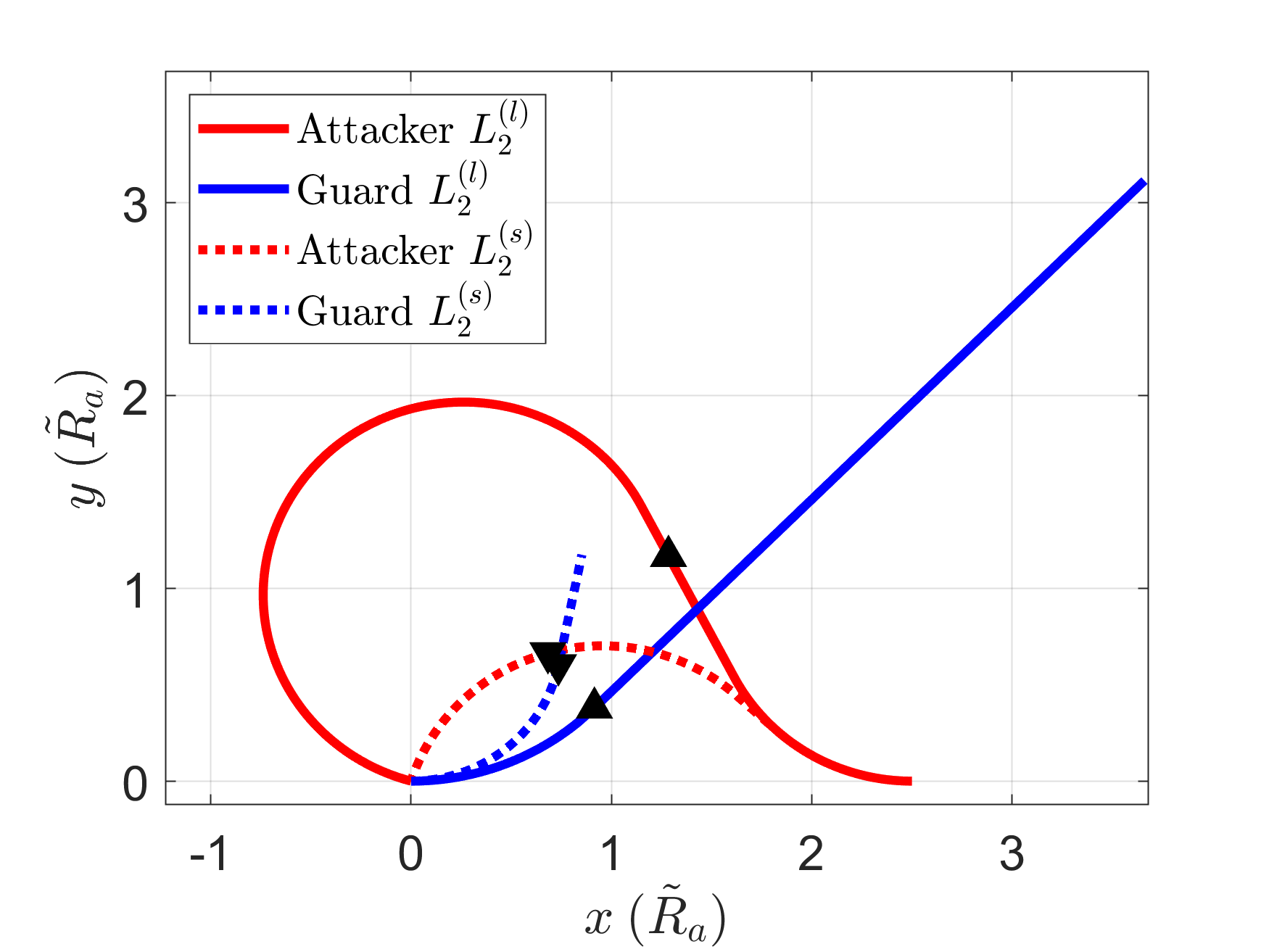}}
\subfigure[$L=\sqrt{8}$.]{\label{fig:bzb_traj2}\includegraphics[width=0.3\textwidth]{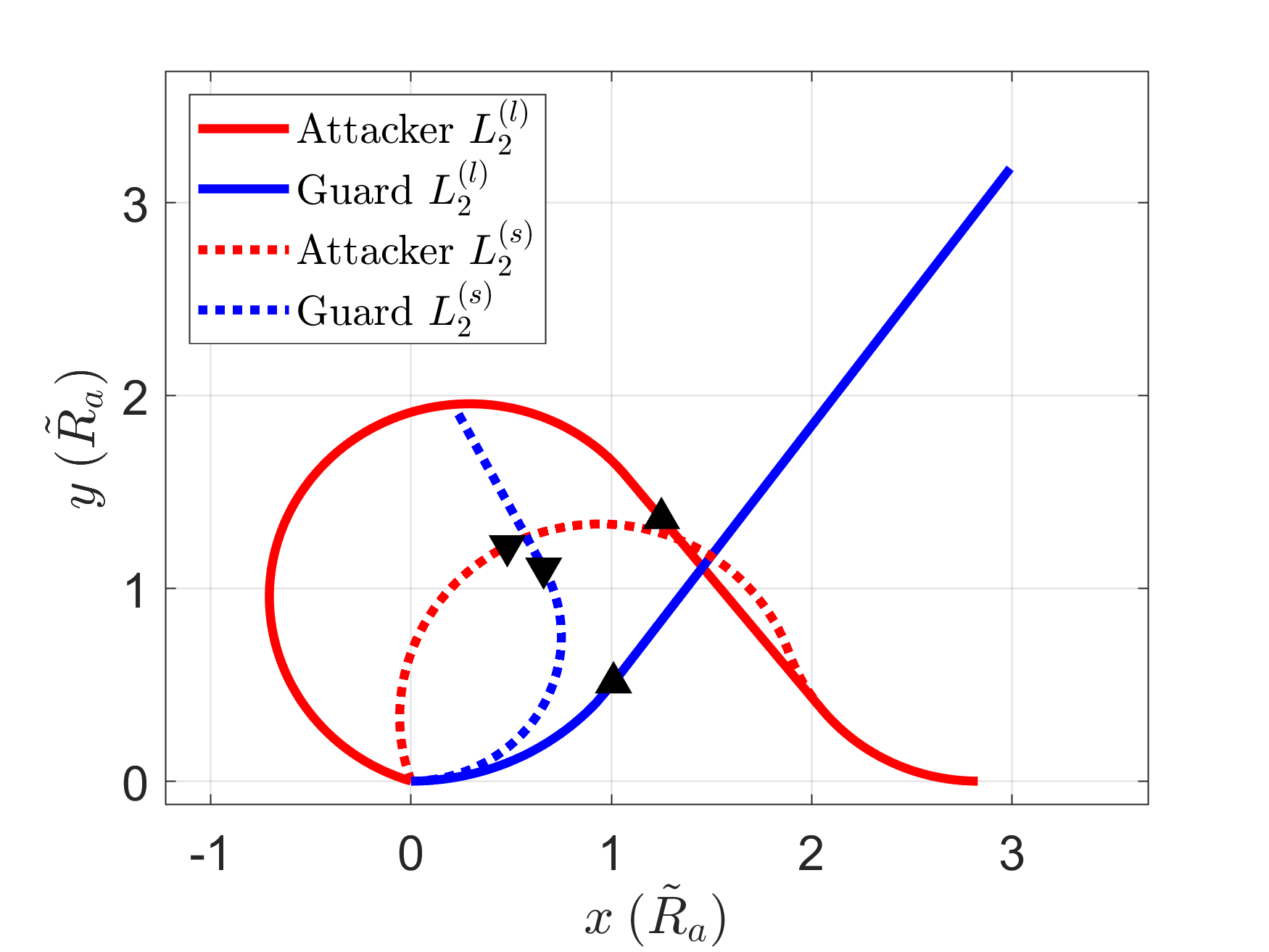}}
\subfigure[$L=3$.]{\label{fig:bzb_traj3}\includegraphics[width=0.3\textwidth]{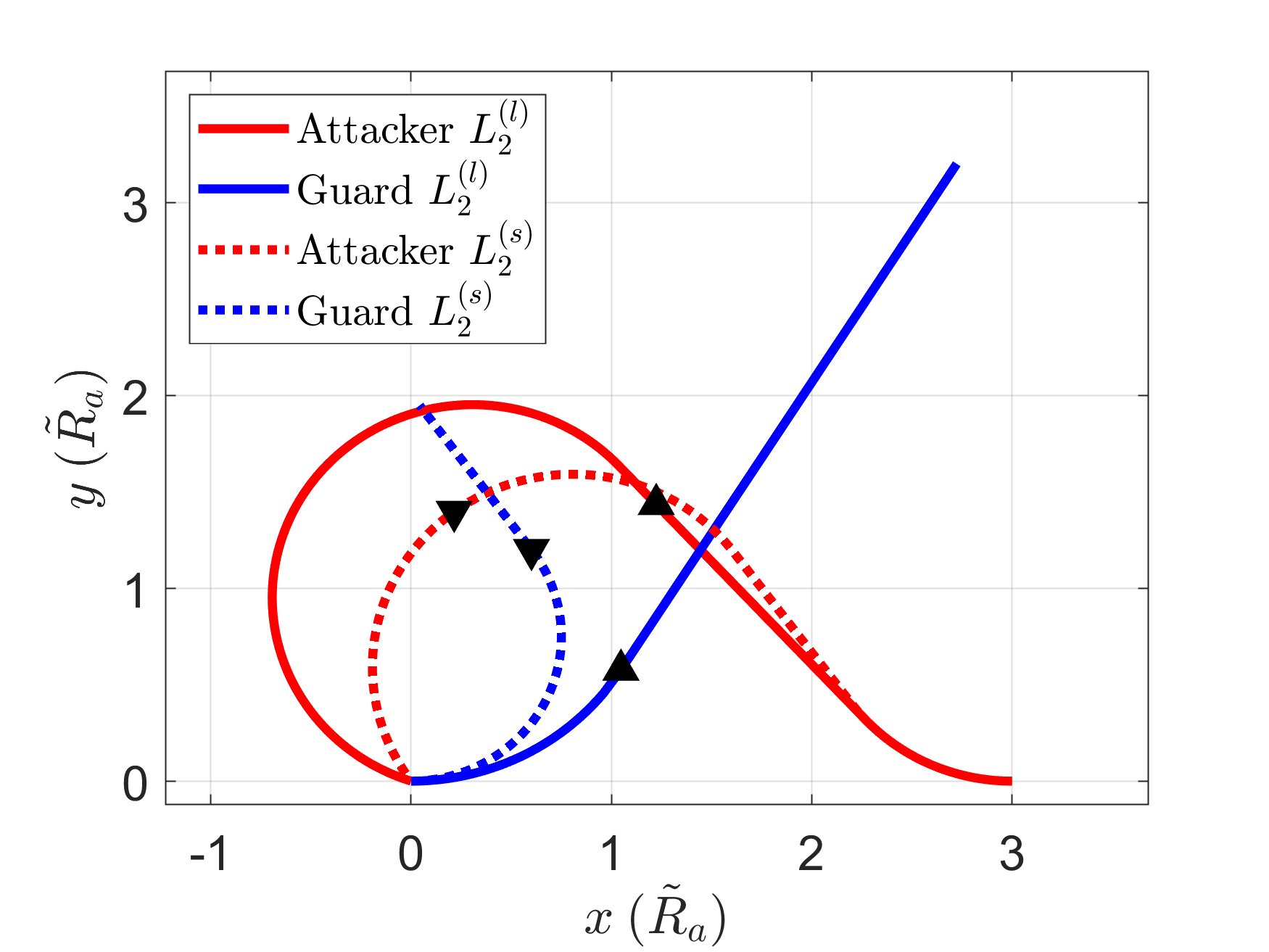}}
\caption{Optimal BZB trajectories.}
\label{fig:bzb_trajectories}
\end{figure}

Figure \ref{fig:bzb_trajectories} shows the trajectories associated with the maxima marked in \figref{fig:miss_vs_L1}. The triangles in this figure mark the players' locations at flyby time.

Returning to our original goal, the above BZB trajectories are useful initial guesses for a numerical CNLP solution to the FBDP problem. To this end, we must supplement them with co-state guess trajectories.  Co-states are the game value sensitivities to variations in state variables that lead to neighboring optimal trajectories. State variations under the restricted BZB optimality are a subspace of the variations allowed in the original problem (for example, moving sideways from the Attacker's final arc is not allowed because he will miss the Target). Co-states computed from BZB trajectories are thus only a projection of the FBDP co-states on this subspace. Nevertheless, we find the following procedure useful:

Let $K(t)=\frac{\partial \pmb{\statevec}_{BZB}(t)}{\partial \pmb{\statevec_0}}$ be the current state to initial state Jacobian, under the BZB constrained optimization \footnote{In writing this, we actually allow a BZB problem with general initial headings. The generalization is straightforward and for the sake of brevity, we avoid it here.}.  Similarly, let $\pmb{\nabla\mu_{BZB}}(t)= \frac{\partial \mu_{BZB}(t)}{\partial \pmb{\statevec}_0}$ be the gradient of the current flyby distance. In particular, $\pmb{\nabla\mu_{BZB}}(t_f)$ is the gradient of the game value w.r.t. the initial state. Note that although $\mu(t)$ is not necessary for the BZB optimization itself, we do need it to approximate FBDP co-states. In analogy with the FBDP definition, we define it as the running minimal players separation (under the BZB constraints):   
\begin{equation}
\mu_{BZB}(t)=\min_{t'\in[0,t]}{||\ri{a}(t')-\ri{g}(t')||}
\end{equation}
With these definitions, we now write the variation of the game value as 

\begin{equation}
\label{eq:costatesapprox1}
\begin{split}
\delta\mu_{BZB}(t_f)\equiv&\pmb{\delta \statevec} ^T\cdot \pmb{\lambda_{kin}}(t)  + \delta \mu_{BZB}\cdot\lambda_{\mu}(t) \\
=&\pmb{\delta \statevec} ^T\cdot \pmb{\lambda_{kin}}(t)  + \lambda_{\mu}(t)\cdot \pmb{\delta \statevec_0 }^T\cdot\pmb{\nabla\mu_{BZB}}(t) \\
=& \pmb{\delta \statevec_0}^T\cdot \pmb{\nabla\mu_{BZB}}(t_f)
\end{split}
\end{equation}

 $\pmb{\lambda_{kin}}$ denotes the vector of kinematic co-states and $\lambda_{\mu}$ is the co-state of the auxiliary flyby distance variable $\mu$. Replacing $\pmb{\delta \statevec_0}$ with ${K}^{-1}(t)\cdot\pmb{\delta \statevec}$ and canceling $\pmb{\delta \statevec}$, we get

\begin{equation}
\label{eq:costatesapprox2}
K(t)^T \cdot\pmb{\lambda_{kin}}(t)  = \pmb{\nabla\mu_{BZB}}(t_f)-\lambda_{\mu}(t)\cdot\pmb{\nabla\mu_{BZB}}(t) 
\end{equation}

For $t\leq t_{fb}$  we use Eq. \ref{eq:approxlmu} to write 
\begin{equation}
\label{eq:costatesapprox3}
K(t)^T \cdot\pmb{\lambda_{kin}}(t)  = \pmb{\nabla\mu_{BZB}}(t_f)-e^{\frac{t-t_{fb}}{\tau}}\cdot\pmb{\nabla\mu_{BZB}}(t) 
\end{equation}
which we solve (using singular value decomposition) for $\pmb{\lambda_{kin}}(t)$ in a least square sense.
For $t>t_{fb}$ , $\lambda_{\mu}(t)\approx1$ and $\mu_{BZB}(t)=\mu_{BZB}(t_f)$, leaving us with  $\pmb{\lambda_{kin}}(t>t_{fb})\approx0$.
This approximation ignores an additive term $\pmb{\lambda_{kin_{null}}}(t)\in Null\left[K(t)^T\right]$ that is lost in the projection from the original to the BZB problem.   

\begin{figure}[bt!]
\centering
\includegraphics[width=0.7\textwidth]{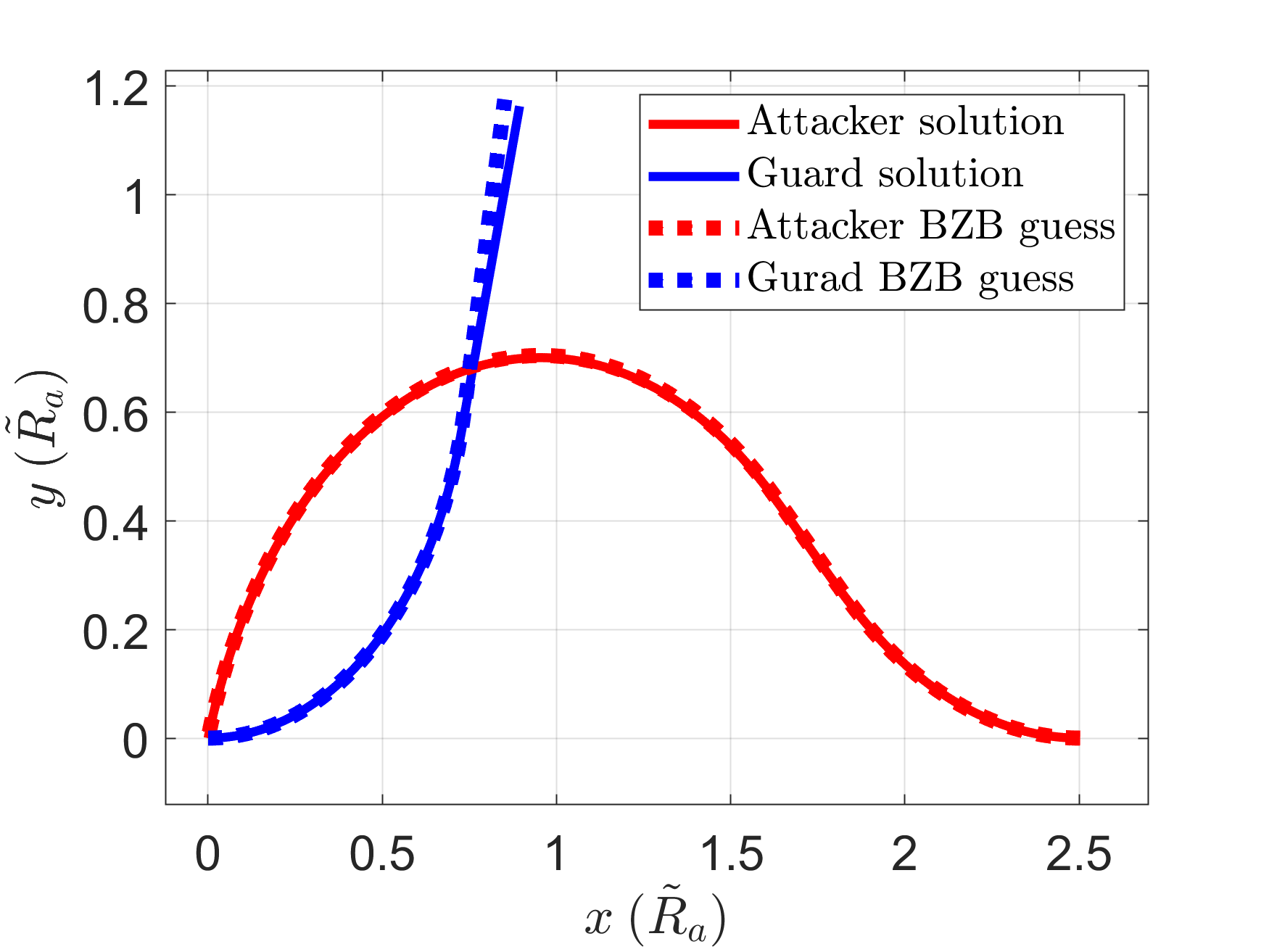}
\caption{BZB guess trajectories comparison with FBDP solution.}
\label{fig:bzb_guess_trajs}
\end{figure}

\begin{figure}[bt!]
\centering    
\subfigure[Guard position co-states.]{\label{fig:bzb_guard_guess_poscostates}\includegraphics[width=0.45\textwidth]{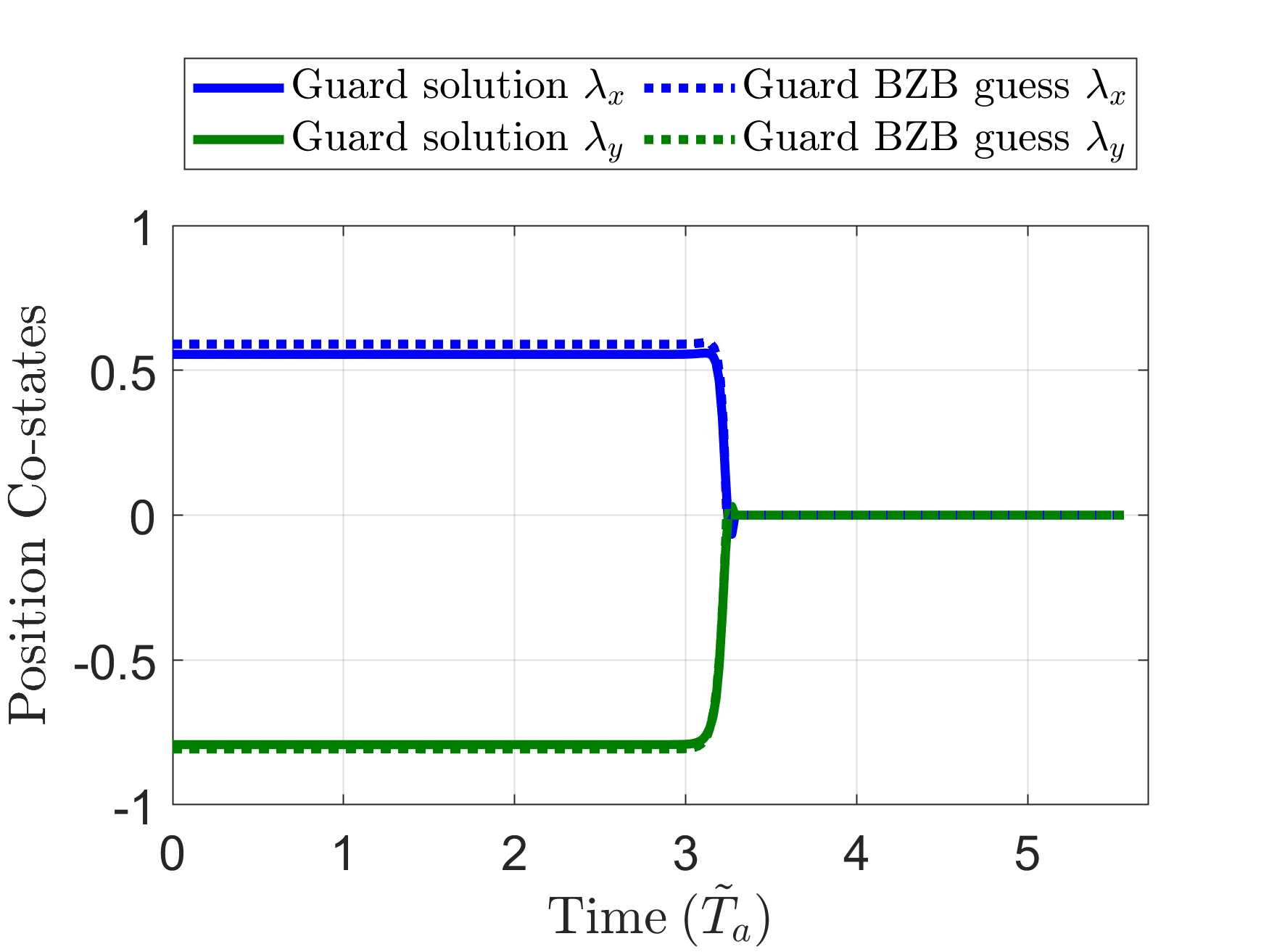}}
\subfigure[Guard velocity co-states.]{\label{fig:bzb_guard_guess_velcostates}\includegraphics[width=0.45\textwidth]{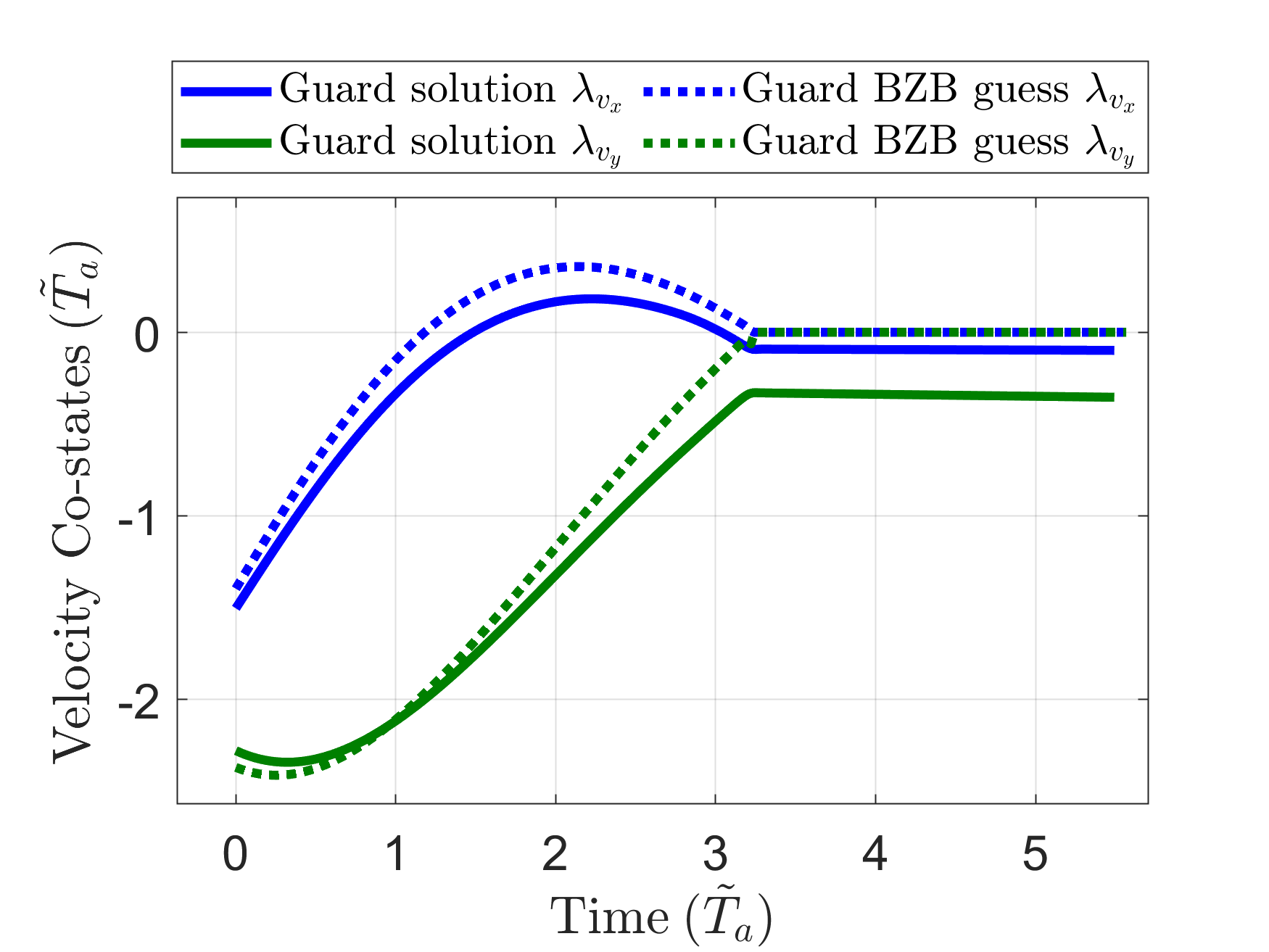}}
\subfigure[Attacker position co-states.]{\label{fig:bzb_attacker_guess_poscostates}\includegraphics[width=0.45\textwidth]{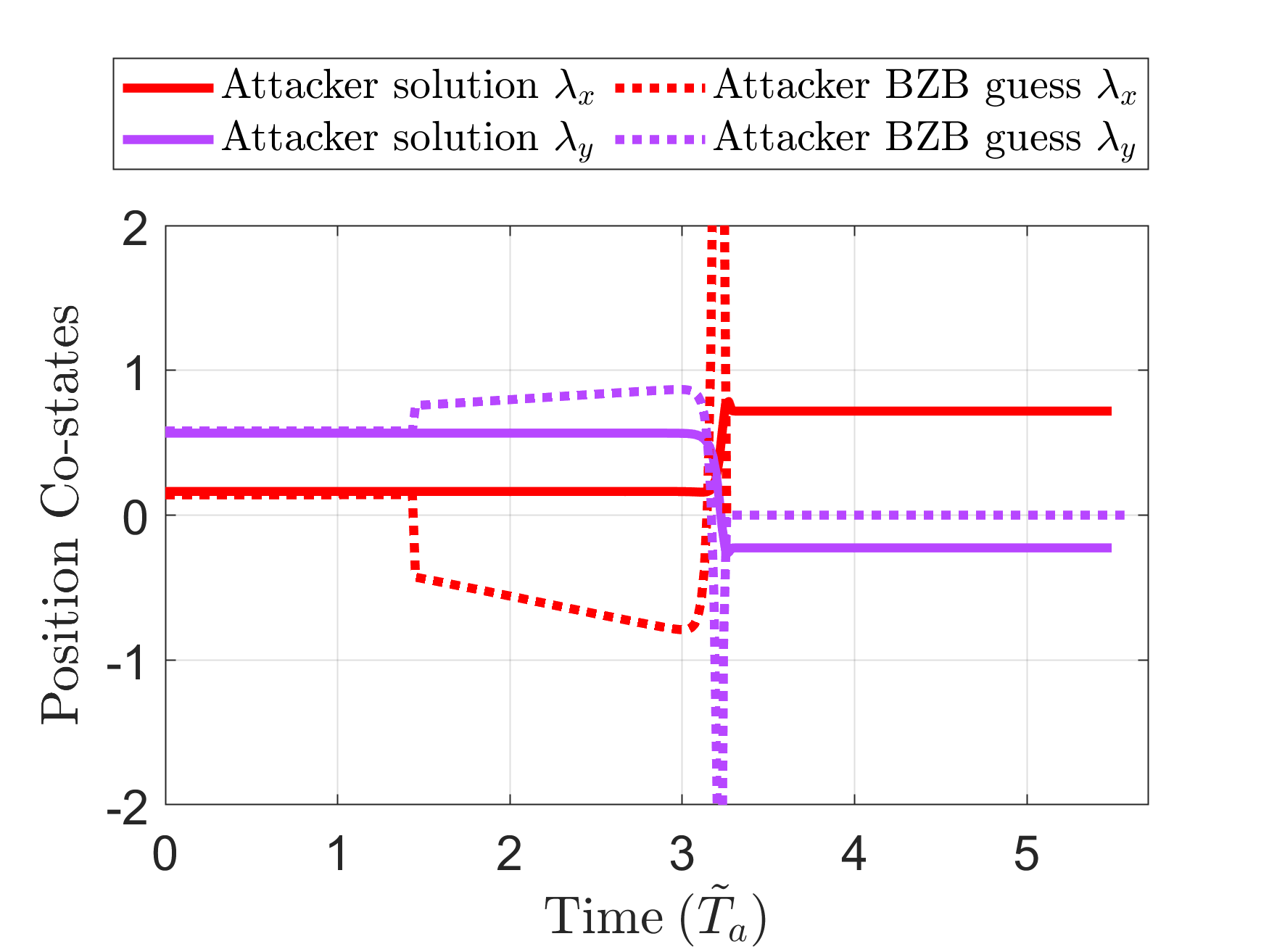}}
\subfigure[Attacker velocity co-states.]{\label{fig:bzb_attacker_guess_velcostates}\includegraphics[width=0.45\textwidth]{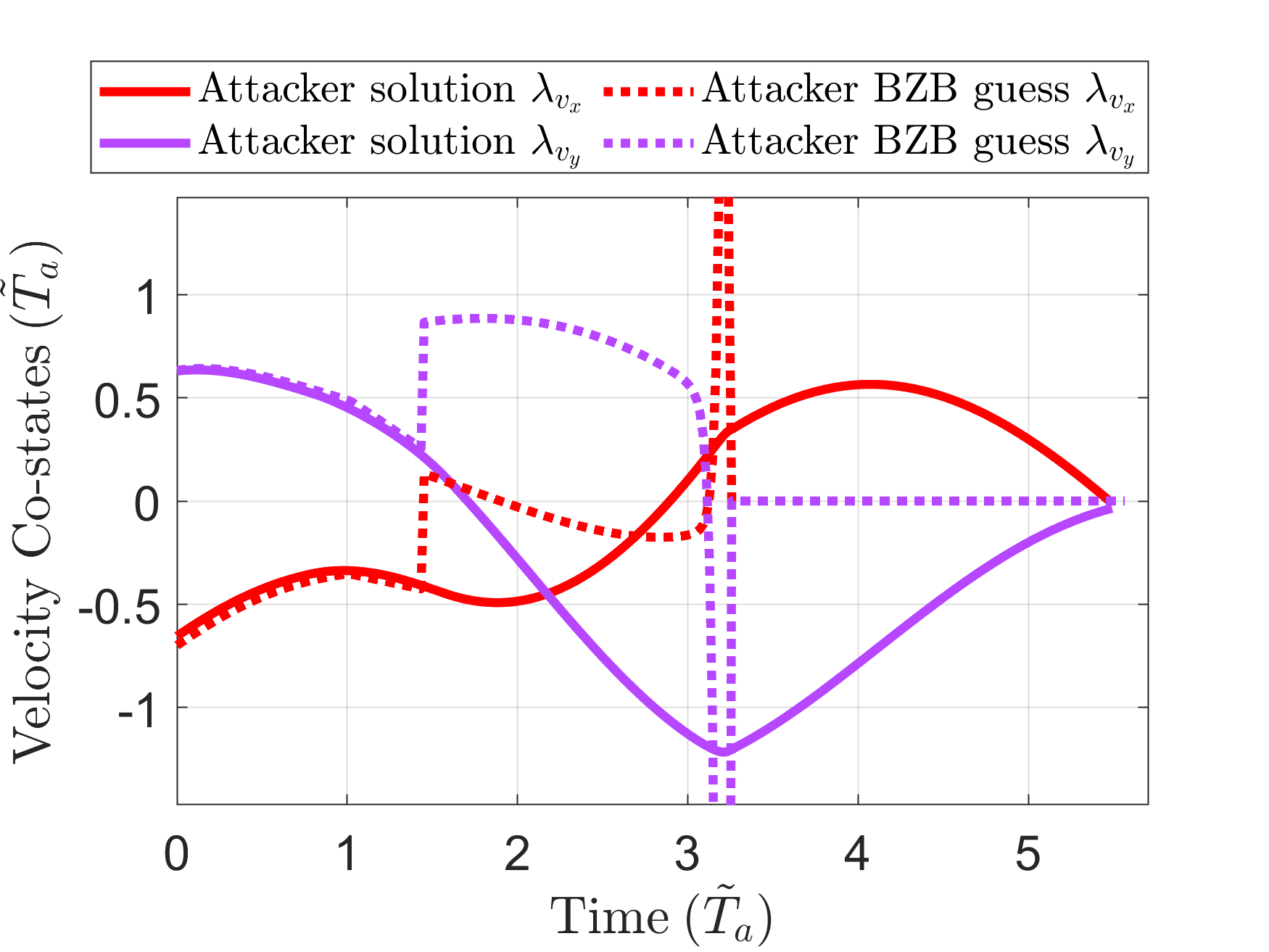}}
\caption{BZB co-states guess comparison with FBDP solution.}
\label{fig:BZB_guess_costates}
\end{figure}

Figure \ref{fig:bzb_guess_trajs} presents a comparison between BZB guess trajectories (dotted curves) and the actual FBDP numerical solution that arrived from it (solid curves). We show the $L=2.5$ local saddle point from \secref{sec:TypeCSolutions} (see \figref{fig:typeC_trajs}). The trajectories agree very well. Figure \ref{fig:BZB_guess_costates} shows a comparison between the guess and the actual co-states. Figures \ref{fig:bzb_guard_guess_poscostates} and \ref{fig:bzb_guard_guess_velcostates} compare the Guard's position and velocity co-states. The guess co-states are a good approximation in this case. Figures \ref{fig:bzb_attacker_guess_poscostates} and \ref{fig:bzb_attacker_guess_velcostates} compare the Attacker's position and velocity co-states. The approximation here is initially good but deteriorates at the end of the first arc and, more so, at the flyby time. This is because the Attacker's subspace of the Jacobian $K(t)$ becomes singular there. Nevertheless, both players' guess co-states satisfy the transversality conditions in Eq. \ref{eq:transversality_bc} (with $\phi_{v_a}=\phi_{v_g}=0$, which is the correct analogy) and are often good enough for CNLP solver convergence. Indeed, the solid curves in \figref{fig:bzb_guess_trajs} resulted from a numerical solution that used the dashed curves as an initial guess.

\bibliography{flyby_distance_pursuit}
\bibliographystyle{plain} 
\end{document}